\documentclass[10pt,twoside,a4paper]{article} 

\usepackage{amsmath}
\usepackage{amsfonts}
\usepackage{cite}

\begin{document}

\begin{center}
{\LARGE \textbf{Extrapolation Methods for Improving the Convergence of 
Oligomer Calculations to the Infinite Chain Limit of 
\textit{Quasi}-Onedimensional Stereoregular Polymers}}
\end{center}

\bigskip

\begin{center}
\textsc{Ernst Joachim Weniger} \\
Institut f{\" u}r Physikalische und Theoretische Chemie \\
Universit{\" a}t Regensburg, D-93040 Regensburg, Germany \\
\texttt{Joachim.Weniger@chemie.uni-regensburg.de} \\ [3\jot]
\textsc{Bernard Kirtman} \\
Department of Chemistry and Biochemistry \\ 
University of California, Santa Barbara, CA 93106-9510, U.S.A. \\
\texttt{kirtman@chem.ucsb.edu}
\end{center}

\begin{center}
{\small 
Submitted to the Special Issue \\ 
``Numerical Methods in Physics, Chemistry and Engineering'' \\ 
of \textit{Computers \& Mathematics with Applications} \\
Guest Editors: T.E.\ Simos, G.\ Avdelas, and J.\ Vigo-Aguiar \\ 
12 April 2000}
\end{center}

\begin{quote}
\textbf{Abstract} --- {\small 
\emph{Quasi}-onedimensional stereoregular polymers as for example 
polyacetylene are currently of considerable interest, both
experimentally as well as theoretically. There are basically two
different approaches for doing electronic structure calculations: One
method, the so-called crystal orbital method, uses periodic boundary
conditions and is essentially based on concepts of solid state
theory. The other method is essentially a quantum chemical method since
it approximates the polymer by oligomers consisting of a finite number
of monomer units, i.e., by molecules of finite size. In this way, the
highly developed technology of quantum chemical molecular programs can
be used. Unfortunately, oligomers of finite size are not necessarily
able to model those features of a polymer which crucially depend of its
in principle infinite extension. However, in such a case extrapolation
techniques can be extremely helpful. For example, one can perform
electronic structure calculations for a sequence of oligomers with an
increasing number of monomer units. In the next step, one then can try
to determine the limit of this sequence for an oligomer of infinite
length with the help of suitable extrapolation methods. Several
different extrapolation methods are discussed which are able to
accomplish an extrapolation of energies and properties of oligomers to
the infinite chain limit. Calculations for the ground state energy of
polyacetylene are presented which demonstrate the practical usefulness
of extrapolation methods.}
\end{quote}

\begin{quote}
\textbf{Keywords} --- {\small
Quantum chemical oligomer calculations, quasi-onedimensional
stereoregular polymers, extrapolation methods, sequence transformations.} 
\end{quote}

\renewcommand{\baselinestretch}{1.15}
\small\normalsize

\setcounter{equation}{0}  
\section{Introduction}
\label{Sec:Intro}

The classic example of a \emph{quasi}-onedimensional stereoregular
polymer is polyacetylene, which can be generated by translating the
repeat unit $-\text{(CH=CH)}-$. Conjugated polymers such as
polyacetylene are currently of great technological interest because of
potential applications in a wide variety of optical and optoelectronic
devices based on their conducting, photophysical, and nonlinear optical
properties \cite{McGeMiMoHee,KirCha1997}. Moreover, the keen interest in 
polyacetylene especially is documented by the fact that there even
exists a monograph with the title \emph{Polyacetylene} \cite{Chi1984}.

Numerous properties involving the electronic structure of these polymers
are not yet understood in a completely satisfactory way. Consequently,
these polymers are not only very interesting from a technological point
of view but are also quite challenging from a basic science point of
view.

There are two approaches to \emph{ab initio} molecular electronic
structure calculations. One may be described as the wavefunction
approach and the other as the density approach (also known as density
functional theory). The starting point for the wavefunction approach is
an effective independent particle treatment of the Schr\"{o}dinger
equation named after Hartree \cite{Har1928} and Fock \cite{Fo1930}. In
the density approach one uses an equation of similar general form due to
Kohn and Sham \cite{KohSha1965} but, in principle, the electron
interactions are fully taken into account. The solution of either
equation is a set of one-electron functions, called orbitals, that can
be used to construct the multi-electron wavefunction (Hartree-Fock) or
the electron density (Kohn-Sham). In either case the orbitals are
expanded in terms of finite linear combination of analytic basis
functions (nowadays usually Gaussians).

The Hartree-Fock approach combined with analytic basis functions leads
to the so-called Roothaan equations \cite{Roo1951}, which can be written
in the form of a generalized matrix eigenvalue problem. For meaningful
results it is necessary for the number of basis functions and, hence,
the dimension of the Roothaan equations to exceed the number of
electrons. This imposes obvious restrictions on the size of the size of
oligomers used to model a polymer, which can be handled by
\emph{ab initio} programs based on the highly developed computational
techniques of molecular quantum chemistry (see for instance
\cite{AnMoMcChClFrLePiVeVr1997,Fri1996}).

If, however, the polymer under consideration is a stereoregular
\emph{quasi}-onedimensional polymer like polyacetylene, it is
nevertheless possible to compute energies or properties on the basis of
quantum mechanics. There are two fundamentally different computational
approaches -- the \emph{crystal orbital} approach and the \emph{finite
oligomer} method -- that have complementary advantages and
disadvantages.

The crystal orbital approach, which is for example discussed in
monographs by Pisani, Dovesi, and Roetti \cite{PisDovRoe1988}, Ladik
\cite{Lad1988}, and Andr\'e, Delhalle, and Br\'edas \cite{AndDelBre1991}
or in review articles by Andr\'e \cite{And1980}, Kert\'esz
\cite{Ker1982}, and Ladik \cite{Lad1999}, treats polymers as solids and 
uses the concepts and techniques of solid state physics such as
translational symmetry and cyclic boundary conditions (see for example
\cite{Kit1971,Kit1963}). This leads to band structure, or crystal 
orbital, calculations whose complexity is limited by the number of
interacting neighboring unit cells. Such calculations have the
undeniable advantage that one only has to consider these interactions
for a representative unit cell since all unit cells are identical. This
is a major simplification over a molecular-type calculation where the
unit cells are translationally inequivalent. Thus, the crystal orbital
approach is in principle well suited for describing those features of a
polymer which crucially depend on its in principle infinite extension.

Nevertheless, there are still some open problems with crystal orbital
calculations. For example, convergence problems in the Hartree-Fock
iterations and the proper handling of long range Coulomb interactions
can be problematic \cite{JacAndCha1999}. Moreover, it is very difficult
to take into account the effect of nonperiodic perturbations such as the
interaction with an electric field. The term in the polarization
response, which is linear in the interacting electric field, determines
the polarizability of the polymer, whereas the nonlinear terms determine
the hyperpolarizabilities. The nonlinear response is connected with
properties such as second harmonic generation, optical rectification,
and the intensity dependent refractive index (see for example
\cite{PraWil1991}) which are of considerable interest in materials
science these days. Despite many years of effort going back to the
1940's and substantial recent progress
\cite{ChaMosAnd1993,ChaOhr1994,ChaMosFriAnd1999,Ott1992,Lad1996,%
OttGuLad1999}, it is fair to say that an incontrovertible crystal
orbital method \cite{JacCha2000,OttGuLad2000} has yet to be presented,
even at a simplified Hartree-Fock level where just the average - rather
than the instantaneous - interelectronic repulsion is taken into
account.

Another problem of the crystal orbital approach is that the explicit
incorporation of electronic correlation is difficult. In the
Hartree-Fock approximation the instantaneous interaction between
electrons is replaced by an orbital-averaged interaction. The error due
to making this approximation is known as the electron correlation error
(see for example \cite{Wil1984,Ful1991,HarMonFre1992}). For some
properties, such as the energy per unit cell, this error is quite small
on a percentage basis. Nonetheless, correlation is very important
because it can easily be the deciding factor between alternative
structural conformations. For other properties, such as
hyperpolarizabilities, the correlation error may be as large as 800 \%
\cite{ToToMel1995}. 

The density method includes electron correlation but the one-electron
potential, which simultaneously takes into account both electron
exchange (Pauli Principle) and correlation, must be approximated. In
practice, the approximations used are local in nature -- they involve
the density, its gradient and, possibly, the Laplacian at each
point. For that reason \cite{vGisSchiGriBaeSniChaKir1999} the method
fails to predict electrical properties (dipole moment, linear
polarizability, hyperpolarizabilites)
\cite{ChaPerGisBaeSniSouRobKir1998,ChaPerJacGisSouRob2000} of spatially
extended systems. One can expect that the same will be true of other
properties of spatially extended systems, such as the mechanical force
constants describing the interaction between different unit cells, but
that has not been checked as yet.

Because of the problems mentioned above, a different approach -- the
\emph{finite oligomer} method -- has gained popularity in the middle to
the late 1980's: Stereoregular \textit{quasi}-onedimensional polymers
are approximated by molecules of finite length
\cite{KiNiPa1983,HuDuCl1988,Kir1988,Cio1988,ViDuWaHuCl1988,WeLi1990,%
CuiKerJia1990,CioWen1992}. In the case of polyacetylene, this means that
the polymer $\text{(CH=CH)}_{\infty}$ is approximated by an oligomer
$\mathrm{H} \! - \! \text{(CH=CH)}_{\mathrm{N}} \! - \! \mathrm{H}$,
where $N$ is a finite integer. The oligomer approach offers many
advantages. The \emph{ab initio} treatment of electron correlation using
any of the techniques available in standard molecular \emph{ab initio}
programs is straightforward. Moreover, the same is true for the
calculation of hyperpolarizabilities, which has yet to be successfully
implemented within the framework of crystal orbital theory even at the
Hartree-Fock level. Molecular-type calculations using finite oligomers
have been done not only at the Hartree-Fock level, but have also
incorporated electron correlation by many-body perturbation theory
\cite{ToToMel1995,ToTodMelHaKi1995,ToToMelRob1995,JacChaAnd1998}. As a 
result, it has become known that the role of electron correlation can be
exceptionally significant for this problem. The same may be said
regarding the importance of the vibrational contribution to the
hyperpolarizability
\cite{ChaLuiDurAndKir2000,ChaKir2000} that arises because the electronic 
response depends upon the nuclear positions. In addition, studies have
been undertaken to determine the effect of interactions in the solid
state with other polymer chains \cite{KirDykCha1999,McWiSo1991} and with
dopants \cite{TaTa1995,Ro-Mo1996,KirSpaCha2000} that either accept or
donate electrons. All of these investigations are readily feasible
because of the availability of quantum chemistry computer codes
developed for ordinary finite molecules.

However, the finite oligomer approach of course also has some
shortcomings which may lead to numerical problems. The basic problem of
the approximation of a \emph{quasi}-linear stereoregular polymer
$\mathrm{A}_{\infty}$ by an oligomer $\mathrm{X} \!- \!
\mathrm{A}_{\mathrm{N}} \!  - \! \mathrm{Y}$, whose terminal free
valences are saturated by suitably terminal groups $\mathrm{X}$ and
$\mathrm{Y}$ (usually hydrogen atoms), is the choice of the number $N$
of repeat units $\mathrm{A}$. Of course, $N$ should be small enough to
permit an accurate quantum chemical calculation of the
oligomer. However, $N$ should also be large enough in order to provide a
sufficiently realistic representation of those features of the polymer
that result from its in principle infinite extension.

Obviously, these two requirements tend to be contradictory which implies
that finite oligomer calculations can suffer badly from slow convergence
to the infinite chain limit. This is where extrapolation methods come
into play: Instead of calculating the intensive energy $E (N)$ or more
generally for an intensive property $P (N)$ of a single oligomer
$\mathrm{X} \! - \! \mathrm{A}_{\mathrm{N}} \!  - \! \mathrm{Y}$, one
can perform calculations for a sequence of oligomers $\mathrm{X} \! - \! 
\mathrm{A}_{\mathrm{N}} \!  - \!
\mathrm{Y}$ with $N = 1, 2, \ldots, N_{\text{max}}$. The sequences 
$E (1)$, $E (2)$, \ldots, $E (N_{\text{max}})$ and $P (1)$, $P (2)$,
\ldots, $P (N_{\text{max}})$ can then be extrapolated to the infinite
chain limits $E ({\infty})$ and $P ({\infty})$, respectively.

Many different approaches for doing the extrapolations are possible. In
this article, we want to use sequence transformations for that
purpose. Numerical techniques for the acceleration of convergence or the
summation of divergent series are as old as calculus. According to Knopp
\cite[p.\ 249]{Kn1964}, the first series transformation was published by
Stirling \cite{St1730} already in 1730, and in 1755 Euler \cite{Eu1755}
published the series transformation which now bears his name. In
rudimentary form, convergence acceleration methods are even older.  
Brezinski \cite[pp.\ 90 - 91]{Br1991a} mentioned that in 1654 Huygens
used a linear extrapolation scheme which is a special case of what we
now call Richardson extrapolation \cite{Ri1927} for obtaining better
approximations to $\pi$, and that in 1674 Seki Kowa, the probably most
famous Japanese mathematician of that period, tried to accomplish this
with the help of the so-called $\Delta^2$ process, which is usually
attributed to Aitken \cite{Ai1926}. Further details on the historical
development of extrapolation methods can be found in an article by 
Brezinski \cite{Bre1996}. 

The modern theory of sequence transformations begins with two articles
by Shanks \cite{Sha1955} and Wynn \cite{Wy1956a}.  Shanks introduced a
sequence transformation which in the case of a power series produces
Pad\'{e} approximants \cite{BaGM1996}, and Wynn showed that this
transformation can be computed conveniently by a nonlinear recursive
scheme called the epsilon algorithm. These two articles appeared at a
time when the first computers became publicly available. The subsequent
revolution created a great demand for powerful numerical tools, and the
two articles by Shanks and Wynn, which offered solutions to practical
problems, had an enormous impact: They stimulated research not only on
Pad\'{e} approximants \cite{BaGM1996} but also on other sequence
transformations, as documented by the large number of monographs
\cite{Bre1977,Bre1978,Bre1980,BreRZa1991,Wi1981,MaSha1983,De1988,%
LiLuSh1995,Wal1996,Bre1997} and review articles \cite{Gu1989,We1989}
which appeared in the last years twenty years or so. Moreover, there is
not only a lot of research on sequence transformations and their
properties, but the number of articles which describe applications also
increases steadily (see for example the list of references in
\cite{We1996a}).

Sequence transformations have already been used for the extrapolation of
oligomer calculations. Initially, they were applied to ground state
energies \cite{WeLi1990,CioWen1992}, but more recently, they have also
been employed for the the computation of static and dynamic linear
polarizabilities \cite{DalOddBis1998}. It is the intention of this
article to provide a sufficiently detailed discussion of those sequence
transformations that are suited for the extrapolation of oligomer
calculations to the infinite limits. As discussed later in more detail,
sequence transformations try to accomplish an improvement of convergence
by detecting and utilizing regularities of the input data. For that
purpose, they have to access the information stored in the later digits
of the input data. However, in the case of oligomer calculations, the  
input data are produced by molecular \emph{ab initio} programs which are
huge packages of FORTRAN code (more than $100~000$ lines of
code). Normally, these programs operate in FORTRAN DOUBLE PRECISION
which corresponds to an accuracy of 14 - 16 decimal digits (depending on
the compiler). Unfortunately, this does not mean that their results have
this accuracy. Numerous internal approximations plus the inevitable
rounding errors reduce the accuracy of the results. This alone is not
necessarily a problem. However, the complexity of the calculations,
which are performed in such a molecular program, makes it impossible to
obtain an estimate of the errors by standard mathematical approaches,
and the only realistic alternative are mathematical experiments. 
Consequently, the emphasis of this manuscript will not be so much on
powerful sequence transformations but rather on sequence transformations
which produce reliable results under the restrictions that are typical
of oligomer calculations.

In Section (\ref{Sec:Eu_Mac}), we give a short discussion of use the
Euler Maclaurin formula for the evaluation of infinite series and show
that sequence transformations also try to eliminate the truncations
errors, albeit by purely numerical means. In Sections (\ref{Sec:AitEps}) 
and (\ref{Sec:RichRho}), we discuss sequence transformations that are
able to acelerate either linear or logarithmic convergence. In
Section (\ref{Sec:Polyac}), we show that oligomer calculations for the
grond state energy of polyacetylene can be extrapolated effectively to
the infinite chain limit by sequence transformations, and finally in
Section (\ref{Sec:SumConclu}) there is a short summary.

\setcounter{equation}{0}
\section{The Euler-Maclaurin Formula and Sequence Transformations}
\label{Sec:Eu_Mac} 

A model example of a slowly convergent series is the Dirichlet series
for the Riemann zeta function:
\begin{equation}
\zeta (z) \; = \; \sum_{\nu = 0}^{\infty} \, (\nu + 1)^{-z} \, .
\label{ZetaSer}
\end{equation}
This series converges provided that $\mathrm{Re} (z) > 1$. However, it
is notorious for extremely slow convergence if $\mathrm{Re} (z)$ is only
slightly larger than 1.

The Euler-Maclaurin formula, which is treated in most books on
asymptotics (see for instance \cite[pp.\ 279 - 295]{Olv1974}), makes it
possible to evaluate $\zeta (z)$ even if the convergence of the
Dirichlet series is extremely slow. Let us consider an infinite series
$\sum_{\nu=0}^{\infty} f (\nu)$, and let us assume that its terms $f
(\nu)$ are \emph{smooth} and \emph{slowly varying} functions of the 
index $\nu$. Then, the integral
\begin{equation}
\int_{M}^{N} \, f (x) \, \mathrm{d} x 
\label{EuMacInt}
\end{equation}
with $M, N \in \mathbb{Z}$ can provide a good approximation to the sum
\begin{equation}
\frac{1}{2} f (M) \, + \, f (M+1) \, + \, \ldots \, + \, f (N-1) \, + \,
\frac{1}{2} f (N) 
\end{equation}
and vice versa. To see this, one only has to approximate the integrand
by step functions.

In the years between 1730 and 1740, Euler and Maclaurin derived
independently correction terms, which ultimately yielded what we now
call the Euler-Maclaurin formula:
\begin{subequations}
\label{EuMaclau}
\begin{align}
\sum_{\nu=M}^{N} \, f (\nu) & \; = \; 
\int_{M}^{N} \, f (x) \, \mathrm{d} x
\, + \, \frac{1}{2} \left[ f (M) + f (N) \right] \notag \\
& \qquad \, + \, \sum_{j=1}^{k} \, \frac {B_{2j}} {(2j)!}
\left[ f^{(2j-1)} (N) - f^{(2j-1)} (M) \right] \, + \, R_k (f) \, , \\
R_k (f) & \; = \; - \, \frac{1}{(2k)!} \,
\int_{M}^{N} \, B_{2k}
\bigl( x - [\mkern - 2.5 mu  [x] \mkern - 2.5 mu ] \bigr) \,
f^{(2k)} (x) \, \mathrm{d} x \, .
\end{align}
\end{subequations}
Here, $[\mkern - 2.5 mu [x] \mkern - 2.5 mu ]$ is the integral part of
$x$, $B_k (x)$ is a Bernoulli polynomial , and $B_m = B_m (0)$ is a
Bernoulli number. 

If we set $M = n+1$ and $N = \infty$, the leading terms of the
Euler-Maclaurin formula provides rapidly convergent approximations to
the truncation error $\sum_{\nu=n+1}^{\infty} f (\nu)$, provided that
$f$ is smooth and a slowly varying function of the index $\nu$. Of
course, it is also necessary that the integral as well as derivatives of
$f$ with respect to $\nu$ can be computed.

The terms $(\nu+1)^{-z}$ of the Dirichlet series (\ref{ZetaSer}) are
obviously smooth functions of the index $\nu$ and they can be
differentiated and integrated easily. Thus, we can apply the
Euler-Maclaurin formula (\ref{EuMaclau}) with $M = n+1$ and $N = \infty$
to the truncation error of the Dirichlet series:
\begin{subequations}
\label{EuMacZeta}
\begin{align}
\sum_{\nu=n+1}^{\infty} \, (\nu+1)^{-z} & \; = \; 
\frac{(n+2)^{1-z}}{z-1} 
\, + \, \frac{1}{2} \, (n+2)^{-z} \notag \\
& \qquad \, + \, 
\sum_{j=1}^{k} \, \frac {(z)_{2j-1} \, B_{2j}} {(2j)!} \,
(n+2)^{-z-2j+1} \, + \, R_k (n,z) \, , \\
R_k (n,z) & \; = \; - \, \frac{(z)_{2k}}{(2k)!} \,
\int_{n+1}^{\infty} \, \frac
{B_{2k} \bigl( x - [\mkern - 2.5 mu  [x] \mkern - 2.5 mu ] \bigr)}
{(1+x)^{z+2k}} \, \mathrm{d} x \, .
\end{align}
\end{subequations}
Here, $(z)_m = z (z+1) \cdots (z+m-1) = \Gamma(z+m)/\Gamma(z)$ is a
Pochhammer symbol.

The expansion (\ref{EuMacZeta}) for the truncation error of the
Dirichlet series (\ref{ZetaSer}) for $\zeta (z)$ is only asymptotic as
$n \to \infty$. Nevertheless, it is numerically extremely useful for
sufficiently large values of $n$, as shown by the following simple
example:
\begin{equation}
\zeta (z) \; \approx \; \sum_{\nu=0}^{n} (\nu+1)^{-z} \, + \, 
\frac{(n+2)^{1-z}}{z-1} \, + \, \frac{1}{2} \, (n+2)^{-z}\, + \, 
\sum_{j=1}^{3} \, \frac {(z)_{2j-1} \, B_{2j}} {(2j)!} \,
(n+2)^{-z-2j+1} \, .
\end{equation}
If we now choose $n = 20$ and $z = 1.01$, we obtain
\begin{subequations}
\begin{align}
\sum_{\nu=0}^{n} (\nu+1)^{-z} & \; = \; 3.599~497~439~829~47 \, , \\
\frac{(n+2)^{1-z}}{z-1} & \; = \; 96.956~241~819~2202 \, , \\
\frac{1}{2 (n+2)^{z}} & \; = \;  
0.220~355~095~043~682 \times 10^{-1} \, , \\
\frac {(z)_{1} \, B_{2}} {2! (n+2)^{z+1}} &
\; = \; 0.168~605~034~844~029 \times 10^{-3} \, , \\
\frac {(z)_{3} \, B_{4}} {4! (n+2)^{z+3}} &
\; = \; - 0.351~266~295~216~895 \times 10^{-7} \, , \\
\frac {(z)_{5} \, B_{6}} {10! (n+2)^{z+5}} & 
\; = \; 0.347~155~401~295~600 \times 10^{-10} \, .
\end{align}
\end{subequations}
The sum of these 6 terms, which is $100.577~943~338~497$, agrees
completely with the ``exact'' result $\zeta (1.01) =
100.577~943~338~497$ obtained by the computer algebra system Maple.

This numerical example shows that a suitable truncation of the
Euler-Maclaurin formula (\ref{EuMaclau}) provides an excellent
approximation to the truncation error $\sum_{\nu=n+1}^{\infty}
(\nu+1)^{-z}$. Here, it should be taken into account that it is
practically impossible to compute $\zeta (1.01)$ with sufficient
accuracy by adding up the terms of the Dirichlet series (\ref{ZetaSer})
successively. We would need $n \approx 10^{600}$ to obtain
$(n+2)^{1-z}/(z-1) = 10^{-4}$ for the leading term on the right-hand
side of Eq.\ (\ref{EuMaclau}), which corresponds to the integral
(\ref{EuMacInt}). Thus, $10^{600}$ terms of the Dirichlet series
(\ref{ZetaSer}) with $z = 1.01$ would lead to an accuracy of only 6
decimal places in the final result.

In view of these excellent results, it looks like an obvious idea to try
to use the Euler-Maclaurin formula also in the case of other slowly
convergent series. Unfortunately, this is not always possible. The
Euler-Maclaurin formula requires that the terms of the series can be
differentiated and integrated with respect to the index. This is only
possible if the terms have a sufficiently simple analytic structure,
which excludes most cases of interest. Moreover, the Euler-Maclaurin
formula cannot be applied in the case of alternating or divergent series
since their terms are neither smooth nor slowly varying functions of the
index.

However, the probably worst drawback of the Euler-Maclaurin formula is
that it is an \emph{analytic} convergence acceleration method. This
means that it cannot be applied if only the numerical values of the
terms of a series are known. This is quite unfortunate since convergence
acceleration techniques are often badly needed in situations in which
apart from a few numerical data very little is known.

These facts give rise to the obvious question whether the
Euler-Maclaurin formula can be generalized in such a way that
approximations to the truncations errors of sequences and series can be
constructed under less restrictive conditions.

Principal tools to overcome convergence problems are \emph{sequence
transformations} which can accomplish something even if no explicit
analytical expression for the truncation error is available. Here, it
should be noted that Pad\'{e} approximants, which have become the
standard tool to overcome convergence problems with power series, are
special sequence transformations since the partial sums of a power
series are transformed into a doubly indexed sequence of rational
functions \cite{BaGM1996}.

Let us assume that we have a sequence $\{ s_n \}_{n=0}^{\infty}$ which
converges to some limit $s = s_{\infty}$. Then, its elements can be
partitioned into the limit $s$ and a remainder $r_n$ according to
\begin{equation}
s_n \; = \; s \, + \, r_n \, , \qquad n \in \mathbb{N}_0 \, .
\end{equation}     
A {\emph{sequence transformation} $\mathcal{T}$ is a rule which
transforms a sequence $\{ s_{n} \}_{n=0}^{\infty}$ into some other
sequence $\{ s'_{n} \}_{n=0}^{\infty}$:
\begin{equation}
\mathcal{T}: \{ s_{n} \}_{n=0}^{\infty} \longmapsto
\{ s'_{n} \}_{n=0}^{\infty} \, .
\end{equation}
A sequence transformation is only useful for our purposes if the
transformed sequence $\{ s'_{n} \}_{n=0}^{\infty}$ converges to the same
limit $s = s_{\infty}$ as the original sequence $\{ s_{n}
\}_{n=0}^{\infty}$. If this is true, then the elements of the
transformed sequence can also be partitioned into the limit $s$ and a
transformed remainder $r'_n$ according to
\begin{equation}
s'_n \; = \; s \, + \, r'_n \, , \qquad n \in \mathbb{N}_0 \, .
\end{equation}
With the exception of some more or less trivial model problems, the
transformed remainders $r'_n$ will be different from zero for all
\emph{finite} values of the index $n$. Thus, it looks as if nothing
substantial is gained by applying a sequence transformation. However,
the new remainders can have much better numerical properties than the
original remainders. 

Accordingly, a sequence transformation is said to \emph{accelerate
convergence} if the transformed remainders $\{ r'_n \}_{n=0}^{\infty}$
vanish more rapidly than the original remainders $\{ r_n
\}_{n=0}^{\infty}$ according to
\begin{equation}
\lim_{n \to \infty} \, \frac{r'_n}{r_n} \; = \;
\lim_{n \to \infty} \, \frac{s'_n - s}{s_n - s} \; = \; 0 \, .
\end{equation}
In the case of some sequence transformations, it is actually possible to
compute the transformed remainders $\{ r'_n \}_{n=0}^{\infty}$
recursively if the numerical values of the original remainders $\{ r_n
\}_{n=0}^{\infty}$ are known \cite{We2000}.

Thus, a sequence transformation tries to construct an approximation to
the actual remainder $r_n$ and to eliminate it from $s_n$. This yields a
new sequence $\{ s'_n \}_{n=0}^{\infty}$ with hopefully superior
numerical properties. In this context, we are confronted with the
practical problem of constructing approximations to the actual
remainders if we do not know explicit analytical expressions for
them. However, even if we know explicit expressions, this normally does
not help. For example, if the input data are the partial sums
\begin{equation}
s_n \; = \; \sum_{k=0}^{n} \, a_k \, ,
\end{equation}
of an infinite series, then the remainders are given by
\begin{equation}
r_n \; = \; - \, \sum_{k=n+1}^{\infty} \, a_k \, .
\label{SerRepRem}
\end{equation}
We automatically have an explicit expression for the remainders if we
have an explicit analytical expression for the terms $a_k$. However,
this does not help. The computation of $r_n$ via the infinite series
(\ref{SerRepRem}) is normally not easier than the direct computation of
the infinite series by adding up the terms successively.  Consequently,
the straightforward computation of the remainders of an infinite series
is normally either not possible or not feasible, and we must use
\emph{indirect} approaches to obtain suitable approximations to the
remainders.

If we want to use a sequence transformation to speed up convergence, the
minimal information, which we need to have, is a \emph{finite} string of
numerical values of input data, for example the sequence elements $s_n$,
$s_{n+1}$, \ldots, $s_{n+k}$ with $k, n \in \mathbb{N}_0$. This alone
is, however, not enough for the construction of an approximation to the
\emph{unknown} truncation error, in particular since there can be an in
principle unlimited variety of different types of remainders. 
Consequently, it should be obvious that a \emph{universal} algorithm for 
the determination and elimination of completely arbitrary unknown 
remainders $\{ r_n \}_{n=0}^{\infty}$ cannot exist. Delahaye and 
Germain-Bonne \cite{DelGerB1980,DelGerB1982} were able to proof this 
rigorously.

Fortunately, in most cases at least some \emph{structural} and
\emph{asymptotic} information about the dependence of the remainders 
$r_n$ on the index $n$ is available, or -- if this is not the case --
some more or less plausible assumptions can be made. We will show in the 
next sections how powerful sequence transformations can be constructed
on the basis of simple assumptions. For example, let us assume that the
elements of a sequence $\{ s_n \}_{n=0}^{\infty}$ can be expressed by
the following series expansion:
\begin{equation}
s_n \; = \; s \, + \, \sum_{j=0}^{\infty} \, c_j \, \varphi_j (n) \, ,
\qquad n \in \mathbb{N}_0 \, .
\label{ModSeqPhi}
\end{equation}
Here, the coefficients $c_j$ are assumed to be unknown, whereas the
functions $\varphi_j (n)$ are assumed to be known, but otherwise they
are essentially arbitrary.

A numerical algorithm can only involve a finite number of arithmetic
operations. Consequently, a \emph{complete numerical} determination of
the remainder $r_n = \sum_{j=0}^{\infty} c_j \varphi_j (n)$ is
impossible since this would require the determination of an infinite
number of unknown coefficients $c_j$. So, the best we can hope for is
the determination of a finite number of unknown coefficients $c_j$. If,
however, the functions $\varphi_j (n)$ decrease with increasing index
$j$, we achieve an acceleration of convergence by eliminating the first
$k$ terms in the infinite series on the right-hand side of
(\ref{ModSeqPhi}), which yields the transformed sequence
\begin{equation}
s'_n \; = \; s \, + \, \sum_{\kappa=0}^{\infty} \, c'_{k+\kappa} \, 
\varphi_{k+\kappa} (n) \, ,
\qquad n \in \mathbb{N}_0 \, .
\label{ModSeqPhiTrans}
\end{equation}
For example, if the functions $\varphi_j (n)$ do not only decrease with
increasing index, but are an \emph{asymptotic sequence} according to
\begin{equation}
\varphi_{j+1} (n) \; = \; \mathrm{o} \bigl( \varphi_j (n) \bigr) \, ,
\qquad n \in \mathbb{}_0 \, ,
\end{equation}
then $s_n$ is according to (\ref{ModSeqPhi}) of order $\mathrm{O} \bigl(
\varphi_0 (n) \bigr)$ as $n \to \infty$, whereas $s'_n$ is according to 
(\ref{ModSeqPhiTrans}) of order $\mathrm{O} \bigl( \varphi_{k} 
(n) \bigr)$, which for sufficiently large values of $k$ leads to a
substantial improvement of convergence.

The elimination of the leading $k$ terms on the right-hand side of
(\ref{ModSeqPhi}) can be accomplished with the help of a sequence
transformation which is exact for the \emph{model sequence}
\begin{equation}
\tilde{s}_n \; = \; \tilde{s} \, + \,
\sum_{j=0}^{k-1} \, \tilde{c}_j \, \varphi_j (n) \, .
\label{Mod_Seq_BH}
\end{equation}
The elements of this model sequence contain $k+1$ unknowns, the limit
$\tilde{s}$ and the $k$ coefficients $\tilde{c}_0$, $\tilde{c}_1$,
$\ldots$, $\tilde{c}_{k-1}$. Since all unknowns occur \emph{linearly},
it follows from Cramer's rule that a sequence transformation can be
constructed which is given as the ratio of two determinants, and which
can determine the limit $\tilde{s}$ in (\ref{Mod_Seq_BH}) if the
numerical values of $k+1$ sequence elements $\tilde{s}_n$,
$\tilde{s}_{n+1}$, $\ldots$, $\tilde{s}_{n+k}$ are available. An
admittedly complicated recursive scheme for this general sequence
transformation, which is determined completely as soon as the functions
$\varphi_j (n)$ are specified, was derived independently by Schneider
\cite{Sch1975}, Brezinski \cite{Bre1980a}, and H{\aa}vie \cite{Haa1979}.

Thus, model sequences $\tilde{s}_n = \tilde{s} + \tilde{r}_n$ with
remainders $\tilde{r}_n$ containing a finite number of terms as in
(\ref{Mod_Seq_BH}) are very useful for the construction of sequence
transformations that are special cases of the general sequence
transformation mentioned above. If the remainders $\{ \tilde{r}_n
\}_{n=0}^{\infty}$ of the model sequence $\{ \tilde{s}_n
\}_{n=0}^{\infty}$ provide sufficiently accurate approximations to the
remainders $\{ r_n \}_{n=0}^{\infty}$ of the sequence 
$\{ s_n \}_{n=0}^{\infty}$ to be transformed, then it can be hoped that 
a sequence transformation, which is exact for the model sequence $\{
\tilde{s}_n \}_{n=0}^{\infty}$, will be able to accelerate the
convergence of $\{ s_n \}_{n=0}^{\infty}$ to its limit $s$. More
detailed discussions of this construction principle as well as many
examples can be found in \cite{BreRZa1991,We1989}.

Before we explicitly construct sequence transformations with the help of
model sequences of the type of (\ref{Mod_Seq_BH}), let us first mention
that in the literature on sequence transformations, simple asymptotic
conditions are used to classify the type of convergence of a
sequence. For example, many practically relevant sequences $\{ s_n
\}_{n=0}^{\infty}$ converging to some limit $s$ can be characterized by
the asymptotic condition
\begin{equation}
\lim_{n \to \infty} \, \frac{s_{n+1} - s}{s_n - s} \; = \; \rho \, ,
\label{LinLogConv}
\end{equation}
which closely resembles the ratio test in the theory of infinite series.
If $\vert \rho \vert < 1$, the sequence is called \emph{linearly}
convergent, and if $\rho = 1$, it is called \emph{logarithmically}
convergent. 

\setcounter{equation}{0}
\section{Aitken Extrapolation and Wynn's Epsilon Algorithm}
\label{Sec:AitEps} 

Let us assume that the remainders of the elements of a sequence $\{ s_n
\}_{n=0}^{\infty}$ consist of a single exponential term:
\begin{equation}
s_n \; = \; s \, + \, c \, \lambda^n \, , \qquad c \ne 0,
\quad \vert \lambda \vert \ne 1 \, , \quad n \in \mathbb{N}_0 \, .
\label{AitModSeq}
\end{equation}
For $n \to \infty$, this sequence obviously converges to its limit $s$
if $0 < \vert \lambda \vert < 1$, and it diverges away from its
generalized limit $s$ if $\vert \lambda \vert > 1$. Thus, if this
sequence converges, it converges linearly according to
(\ref{LinLogConv}). Moreover, this model sequence contains the partial
sums $\sum_{\nu=0}^{n} (-z)^{\nu} = [1-(-z)^{n+1}]/[1+z]$ of the
geometric series as special cases.

If we consider $s$, $c$, and $\lambda$ in (\ref{AitModSeq}) as unknowns
of the linear system $s_{n+j} = s + c \lambda^{n+j}$ with $j = 0, 1, 2$,
then we can easily construct a sequence transformation, which is able to
determine the (generalized) limit $s$ of the model sequence
(\ref{AitModSeq}) from the numerical values of three consecutive
sequence elements $s_n$, $s_{n+1}$ and $s_{n+2}$. A short calculation
shows that
\begin{equation}
\mathcal{A}_{1}^{(n)} \; = \; s_n \, - \,
\frac{[\Delta s_n]^2}{\Delta^2 s_n} \, ,  \qquad n \in \mathbb{N}_0 \, ,
\label{AitFor_1}
\end{equation}
is able to determine the limit $s$ of the model sequence
(\ref{AitModSeq}) according to $\mathcal{A}_{1}^{(n)} = s$. The forward
difference operator $\Delta$ in (\ref{AitFor_1}) is defined by its
action on a function $g = g (n)$:
\begin{equation}
\Delta g (n) \; = \; g (n+1) \, - \, g (n) \, .
\end{equation}           

The $\Delta^2$ formula (\ref{AitFor_1}) is certainly one of the oldest
sequence transformations. It is usually attributed to Aitken
\cite{Ai1926}, but it is actually much older. As already mentioned in
Section \ref{Sec:Intro}, it was used in 1674 by Seki Kowa, the probably
most famous Japanese mathematician of that period \cite[pp.\ 90 -
91]{Br1991a}, and according to Todd \cite[p.\ 5]{Tod1962} it was in principle
already known to Kummer \cite{Kum1837}.  

The power and practical usefulness of Aitken's $\Delta^2$ formula is of
course limited since it is designed to eliminate only a single
exponential term from the elements of the model sequence
(\ref{AitModSeq}). However, the quantities ${\cal A}_1^{(n)}$ can again
be used as input data in (\ref{AitFor_1}). Hence, the $\Delta^2$ process
can be iterated, yielding the following nonlinear recursive scheme
\cite[Eq.\ (5.1-15)]{We1989}:
\begin{subequations}
\label{It_Aitken}
\begin{eqnarray}
\lefteqn{} \nonumber \\
{\cal A}_0^{(n)} & = & s_n \, , \qquad n \in \mathbb{N}_0 \, , \\
{\cal A}_{k+1}^{(n)} & = & {\cal A}_{k}^{(n)} -
\frac
{\bigl[\Delta {\cal A}_{k}^{(n)}\bigr]^2}
{\Delta^2 {\cal A}_{k}^{(n)}} \, , \qquad k, n \in \mathbb{N}_0 \, .
\end{eqnarray}
\end{subequations}%
In the case of doubly indexed quantities like $\mathcal{A}_{k}^{(n)}$,
it will always be assumed that the difference operator $\Delta$ only
acts on the superscript $n$ but not on the subscript $k$:
\begin{equation}
\Delta \mathcal{A}_{k}^{(n)} \; = \;
\mathcal{A}_{k}^{(n+1)} \, - \, \mathcal{A}_{k}^{(n)} \, .
\end{equation}  
The numerical performance of Aitken's iterated $\Delta^2$ process was
studied in \cite{We1989,SmiFor1982}. Concerning the theoretical
properties of Aitken's iterated $\Delta^2$ process, very little seems to
be known. Hillion \cite{Hil1975} was able to find a model sequence for
which the iterated $\Delta^2$ process is exact. He also derived a
determinantal representation for $\mathcal{A}_k^{(n)}$. However,
Hillion's expressions for $\mathcal{A}_k^{(n)}$ contain explicitly the
lower order transforms $\mathcal{A}_0^{(n)}, \ldots,
\mathcal{A}_{k-1}^{(n)}, \ldots, \mathcal{A}_0^{(n+k)}, \ldots,
\mathcal{A}_{k-1}^{(n+k)}$. Consequently, it seems that Hillion's result 
\cite{Hil1975} -- although interesting from a formal point of view -- 
cannot help much to understand the theoretical properties of
$\mathcal{A}_{k}^{(n)}$.

A more detailed discussion of Aitken's iterated $\Delta^2$ process as
well as additional references can for instance be found in Section 5 of
\cite{We1989} or in \cite{We2000}. The iteration of other sequence
transformations is discussed in \cite{We1991}.

An obvious generalization of the model sequence (\ref{AitModSeq}) would
be the following model sequence which contains $k$ exponential terms:
\begin{equation}
s_n \; = \; s \, + \, \sum_{j=0}^{k-1} \, c_j \, \lambda_{j}^{n} \, ,
\qquad n \in \mathbb{N}_0 \, .
\label{EpsModSeq1}
\end{equation}
Here, we assume that $\vert \lambda_{0} \vert > \vert \lambda_{1} \vert
> \ldots > \vert \lambda_{k-1} \vert$ holds. Although the $\Delta^2$
process (\ref{AitFor_1}) is by construction exact for the model sequence
(\ref{AitModSeq}), the iterated $\Delta^2$ process (\ref{It_Aitken}) is
not exact for the model sequence (\ref{EpsModSeq1}). Instead, this is
true for Wynn's epsilon algorithm \cite{Wy1956a}, which is another
generalization of the $\Delta^2$ process (\ref{AitFor_1}) and which
corresponds to the following nonlinear recursive scheme:
\begin{subequations}
\label{eps_al}
\begin{eqnarray}
\epsilon_{-1}^{(n)} & \; = \; & 0 \, ,
\qquad \epsilon_0^{(n)} \, = \, s_n \, ,
\qquad  n \in \mathbb{N}_0 \, , \\
\epsilon_{k+1}^{(n)} & \; = \; & \epsilon_{k-1}^{(n+1)} \, + \,
1 / [\epsilon_{k}^{(n+1)} - \epsilon_{k}^{(n)} ] \, ,
\qquad k, n \in \mathbb{N}_0 \, .
\end{eqnarray}
\end{subequations}
Only the elements $\epsilon_{2k}^{(n)}$ with \emph{even} subscripts
provide approximations to the limit $s$ of the sequence $\{ s_n
\}_{n=0}^{\infty}$ to be transformed. The elements
$\epsilon_{2k+1}^{(n)}$ with \emph{odd} subscripts are only auxiliary
quantities which diverge if the whole process converges.

A straightforward calculation shows that $\mathcal{A}_{1}^{(n)} =
\epsilon_{2}^{(n)}$. Accordingly, Aitken's iterated $\Delta^2$ process
may also be viewed as an iteration of $\epsilon_{2}^{(n)}$. However, for
$k > 1$, $\mathcal{A}_{k}^{(n)}$ and $\epsilon_{2k}^{(n)}$ are in
general different, although they of course have similar properties in
convergence acceleration and summation processes.

If the input data $s_n$ are the partial sums $f_n (z) = \sum_{\nu=0}^{n}
\gamma_{\nu} z^{\nu}$ of the (formal) power series for some function 
$f (z)$, $s_n = f_n (z)$, then Wynn \cite{Wy1956a} could show that his
epsilon algorithm produces Pad\'{e} approximants according to
\begin{equation}
\epsilon_{2 k}^{(n)} \; = \; [ n + k / k ]_f (z) \, .
\label{Eps_Pade}
\end{equation}
For Pad\'e approximants, we use the notation of Baker and Graves-Morris
\cite{BaGM1996}, i.e., the Pad\'e approximant $[l/m]_f (z)$ to some
function $f (z)$ is the ratio of two polynomials $P_l (z) = p_0 + p_1 z
+ \ldots + p_l z^l$ and $Q_m (z) = 1 + q_1 z + \ldots + q_m z^m$ of
degrees $l$ and $m$ in $z$ according to
\begin{equation}
[l/m]_f (z) \; = \; P_l (z) / Q_m (z) \, .
\label{Def_PA}
\end{equation}      

Since the epsilon algorithm can be used for the computation of Pad\'e
approximants according to (\ref{Eps_Pade}), it is discussed in books on
Pad\'e approximants such as \cite{BaGM1996}. However, there is also an
extensive literature which deals directly with the epsilon algorithm. On
p.\ 120 of Wimps book \cite{Wi1981} it is mentioned that over 50
articles on the epsilon algorithm were published by Wynn alone, and at
least 30 articles by Brezinski. As a fairly complete source of
references Wimp recommends Brezinski's first book
\cite{Bre1977}. However, this book was published in 1977, and since then
many more articles on the epsilon algorithm have been
published. Moreover, the epsilon algorithm is not restricted to scalar
sequences but can be generalized to cover for example vector
sequences. A very recent review of these developments can be found in
\cite{GMRobSal2000}.

Aitken's iterated $\Delta^2$ process (\ref{It_Aitken}) as well as Wynn's
epsilon algorithm (\ref{eps_al}) are powerful accelerators for sequences
which according to (\ref{LinLogConv}) converge linearly. However,
contrary to a widespread misconception, the epsilon algorithm is not
necessarily the most efficient convergence accelerator for linear
convergence. In 1973 Levin \cite{Lev1973} introduced a sequence
transformation which uses as input data not only a finite substring
$s_{n}$, $s_{n+1}$, \ldots, $s_{n+k}$ of the sequence to be transformed,
but also explicit estimates $\omega_{n}$, $\omega_{n+1}$, \ldots,
$\omega_{n+k}$ for the corresponding remainders $r_{n}$, $r_{n+1}$,
\ldots, $r_{n+k}$. The explicit incorporation of the information
contained in the remainder estimates into the convergence acceleration
or summation process makes Levin's sequence transformation
\cite{Lev1973} as well as closely related generalizations
\cite{We1989,Hom2000} potentially very powerful. Some recent
applications of Levin-type transformations plus additional references
can for example be found in
\cite{We1996a,We1996b,We1996c,We1996d,JeMoSoWe1999,JenBecEWenSof2000}.

However, the use of explicit remainder estimates is not only the major
advantage but also the major weakness of Levin-type transformations. If
it is possible to find remainder estimates that are good approximations
to the actual remainders of the sequence to be transformed, then
experience indicates that Levin-type transformations can be extremely
powerful. In contrast, Levin-type transformations will perform quite
poorly if good remainder estimates cannot be found. 

If we want to extrapolate oligomer calculations, then we have to
struggle not only with potentially large numerical inaccuracies in the
input data but also with the fact that oligomers with a small number of
repeat units may behave very irregularly. As a result large
\emph{prediction} errors can occur leading to poor performance of
Levin-type transformations.

Wynn's epsilon algorithm (\ref{eps_al}) only uses the elements of the
sequence to be transformed as input data and it is not possible to feed
in any additional information into the transformation
process. Consequently, the epsilon algorithm is not as powerful as
Levin-type transformation can be under optimal conditions. On the other
hand, there is no prediction error. Moreover, as for example discussed
in \cite[Section~15.2]{We1989}, the epsilon algorithm is also remarkably
insensitive to rounding errors and it can tolerate input data which
either have a low relative accuracy or which behave in a comparatively
irregular way. Thus, the real strength of the epsilon algorithm is its
robustness. Due to its remarkable robustness, the epsilon algorithm is
often able to produce meaningful and reliable results in situations in
which other sequence transformations, which are in principle more
powerful, fail. Experience indicates that Aitken's iterated $\Delta^2$
process (\ref{It_Aitken}) is also less robust and more susceptible to
rounding errors than Wynn's epsilon algorithm.

\setcounter{equation}{0}
\section{Richardson Extrapolation, Wynn's Rho Algorithm, and Related
matters} 
\label{Sec:RichRho}

In Section \ref{Sec:AitEps}, we discussed sequence transformations which
work in the case of sequences that converge \emph{linearly} according to
(\ref{LinLogConv}).  Typical examples of linearly convergent sequences
are the partial sums of power series with a nonzero but finite radius of
convergence. The computational problems, which occur in that case, can
be studied via the Gaussian hypergeometric series
\begin{equation}
{}_2 F_1 (a, b; c; z) \; = \; \sum_{m=0}^{\infty} \,
\frac {(a)_m (b)_m} {(c)_m m!} \, z^m \, .
\label{Ser_2f1}
\end{equation}
This series terminates if either $a$ or $b$ is a negative integer.
Otherwise it converges linearly in the interior of the unit circle,
i.e., for $\vert z \vert < 1$, and it diverges for $\vert z \vert >
1$. Moreover, the convergence of this series can for most values of the
parameters $a, b, c \in \mathbb{C}$ be accelerated quite efficiently by
sequence transformations \cite{We2000a}.

Generally speaking, the acceleration of linear convergence is
comparatively simple, both theoretically and practically, as long as
$\rho$ in (\ref{LinLogConv}) is not too close to 1 (this would be a
borderline case to logarithmic convergence). With the help of
Germain-Bonne's formal theory of convergence acceleration
\cite{GerB1973} and its extension \cite[Section 12]{We1989}, it can be
decided rigorously whether a sequence transformation is capable of
accelerating linear convergence or not. Moreover, many sequence
transformations are capable of accelerating linear convergence
effectively.

The situation is much more difficult in the case of logarithmic
convergence. It follows from the expansion (\ref{EuMacZeta}) for the
truncation error that the partial sums $\sum_{\nu=0}^{n} (\nu+1)^{-z}$
of the Dirichlet series (\ref{ZetaSer}) converge logarithmically to
$\zeta (z)$. As already discussed in Section \ref{Sec:Eu_Mac}, the
convergence of the Dirichlet series can become so slow that the
evaluation of $\zeta (z)$ by successively adding up the terms
$(\nu+1)^{-z}$ is practically impossible. Of course, analogous problems
can also occur in the case of other logarithmically convergent sequences
and series,

Moreover, there are also some principal theoretical problems. Delahaye
and Germain-Bonne \cite{DelGerB1980,DelGerB1982} showed that no sequence
transformation can exist which is able to accelerate the convergence of
\emph{all} logarithmically convergent sequences. Consequently, in the
case of logarithmic convergence the success of a convergence
acceleration process cannot be guaranteed unless additional information
is available. Also, an analogue of Germain-Bonne's beautiful formal
theory of the acceleration of linear convergence \cite{GerB1973} and its
extension \cite[Section 12]{We1989} cannot exist.

In spite of these complications, many sequence transformations are known
which work reasonably well at least for suitably restricted subsets of
the class of logarithmically convergent sequences. Examples, which will
be discussed later, are Richardson extrapolation \cite{Ri1927}, Wynn's
rho algorithm \cite{Wy1956b} and its iteration \cite[Section 6]{We1989},
as well as Osada's modification of the rho algorithm \cite{Os1990}, and
the modification of the $\Delta^2$ process by Bj{\o}rstad, Dahlquist,
and Grosse \cite{BjDaGr1981}. Nevertheless, there is a considerable
amount of theoretical and empirical evidence that sequence
transformations speed up logarithmic convergence in general less
efficiently than linear convergence.

Another problem with logarithmic convergence is that numerical
instabilities are more likely than in the case of linear convergence. A
sequence transformation can only accelerate convergence if it succeeds
in extracting some information about the index-dependence of the
truncation errors from a finite set $s_n$, $s_{n+1}$, $\ldots$,
$s_{n+k}$ of input data. Normally, this is done by forming arithmetic
expressions involving higher weighted differences. However, forming
higher weighted differences is a potentially unstable process which can
easily lead to a serious loss of significant digits or even to
completely nonsensical results. If the input data are the partial sums
of a \emph{strictly alternating} series, the formation of higher
weighted differences is normally a remarkably stable process. If,
however, the input data all have the \emph{same sign}, numerical
instabilities due to cancellation are quite likely, in particular if
convergence is very slow. Thus, if the sequence to be transformed
converges logarithmically, numerical instabilities are a serious problem
and at least some loss of significant digits is to be expected.

In \cite{JeMoSoWe1999}, a combined nonlinear-condensation transformation
was described which first converts a slowly convergent monotone series
into a strictly alternating series. Since the subsequent transformation
of a strictly alternating series is a remarkable stable process, it is
in this way possible to evaluate special functions, that are defined by
extremely slowly convergent series, not only relatively efficiently but
also very close to machine accuracy. Unfortunately, an analogous
approach is not possible if oligomer calculations produce
logarithmically convergent sequences: The conversion of the monotone
series to the alternating series requires that series terms with large
indices can be computed. This is of course not possible in the case of
oligomer calculations.

For the construction of sequence transformations, the standard
interpolation and extrapolation methods of numerical mathematics are
quite helpful. Let us assume that the values of a function $f(x)$ are
only known at some discrete points $x_0 < x_1 < \cdots < x_m$, and that
we want to estimate the value of $f$ at some point $\xi \notin \{ x_0,
x_1, \ldots , x_m \}$. If $x_0 < \xi < x_m$, this problem is called 
\textit{interpolation}, and if either $\xi < x_0$ or $x_m < \xi$, this
problem is called \textit{extrapolation}. These topics are discussed in any
book on numerical analysis. More specialized treatments
can for instance be found in \cite{Dav1975,Joy1971}.

For the construction of extrapolation processes, we postulate the
existence of a function $\mathcal{S}$ of a continuous variable which
coincides on a set of discrete arguments $\{ x_n \}_{n=0}^{\infty}$
with the elements of the sequence $\{ s_n \}_{n=0}^{\infty}$ to be
transformed: 
\begin{equation}
\mathcal{S} (x_n) \; = \; s_n \, , \qquad n \in \mathbb{N}_0 \, .
\end{equation}
This ansatz reduces the problem of accelerating the convergence of a
sequence to an extrapolation problem. If a finite string $s_m, s_{m+1},
\ldots , s_{m+k}$ of $k+1$ sequence elements is known one can construct
an approximation $\mathcal{S}_k (x)$ to $\mathcal{S} (x)$ which
satisfies the $k+1$ interpolation conditions
\begin{equation}
\mathcal{S}_k (x_{m+j}) \; = \; s_{m+j} \, ,
\qquad 0 \le j \le k \, .  
\label{IntPolCond}
\end{equation}
In the next step, the value of $\mathcal{S}_k (x)$ has to be determined
for $x \to x_{\infty}$. If this can be done and if there exists a
function $\mathcal{S} (x)$ which can be approximated at least locally by
a suitable set of interpolating functions, one can expect that
$\mathcal{S}_k (x_{\infty})$ will provide a better approximation to the
limit $s = s_{\infty}$ of the sequence $\{ s_n \}_{n=0}^{\infty}$ than
the last sequence element $s_{m+k}$ which was used for the construction
of $\mathcal{S}_k (x)$.

The most common interpolating functions are either polynomials or
rational functions. These two sets will also lead to different
convergence acceleration methods. If a sequence transformation is
constructed on the basis of polynomial interpolation, it is implicitly
assumed that the $k$-th order approximant ${\cal S}_k (x)$ is a
polynomial of degree $k$ in $x$:
\begin{equation}
\mathcal{S}_k (x) \; = \, c_0 \, + \, c_1 x \, + \, \cdots 
\, + \, c_k x^k \, . 
\label{InterpolPol}
\end{equation}
For polynomials, the most natural extrapolation point is $x = 0$.
Consequently, the interpolation points $x_n$ have to satisfy the
conditions 
\begin{subequations}
\label{x_n2zero}
\begin{eqnarray}
& x_0 > x_1 > \cdots > x_m > x_{m+1} > \cdots > 0 \, ,  
\\
& {\displaystyle \lim_{n \to \infty} \; x_n \; = \; 0} \, . 
\end{eqnarray}
\end{subequations}
The choice $x_{\infty}=0$ as the extrapolation point implies that the
approximation to the limit is to be identified with the constant term
$c_0$ of the polynomial (\ref{InterpolPol}).

Several different methods for the construction of interpolating
polynomials ${\cal S}_k (x)$ are known in the mathematical
literature. Since only the constant term of a polynomial $\mathcal{S}_k
(x)$ has to be computed and since in most applications it is desirable
to compute simultaneously a whole string of approximants $\mathcal{S}_0
(0), \mathcal{S}_1 (0), \mathcal{S}_2 (0), \ldots$, the most economical
choice is Neville's scheme \cite{Nev1934} for the recursive computation
of interpolating polynomials. If we set $x=0$, Neville's algorithm
reduces to the following 2-dimensional linear recursive scheme \cite[p.\
6]{Bre1978}:
\begin{subequations}
\label{RichAl}
\begin{eqnarray}
\mathcal{N}_0^{(n)} (s_n , x_n ) & = \; & s_n \, , 
\quad n \in \mathbb{N}_0 \, , 
\\
\mathcal{N}_{k+1}^{(n)} (s_n , x_n ) & = &
\frac
{x_n \mathcal{N}_{k}^{(n+1)} (s_{n+1} , x_{n+1} ) \, - \,
x_{n+k+1} \mathcal{N}_{k}^{(n)} (s_n , x_n ) }
{x_n \, - \, x_{n+k+1}} \; ,
\quad k,n \in \mathbb{N}_0 \, . 
\end{eqnarray}
\end{subequations}
In the literature on convergence acceleration this variant of Neville's
recursive scheme is usually called Richardson extrapolation
\cite{Ri1927}. Obviously, $\mathcal{N}_k^{(n)} (s_n , x_n)$ is exact if 
elements $s_n$ of the sequence to be transformed are polynomials of
degree $k$ in the interpolation points $x_n$:
\begin{equation}
s_n \; = \; s \, + \, \sum_{j=0}^{k-1} \; c_j x_n^{j+1} \, ,
\qquad k,n \in \mathbb{N}_0 \, .
\end{equation}

As is well known, some functions can be approximated by polynomials only
quite poorly, but by rational functions they can be approximated very
well. Consequently, it is to be expected that at least for some
sequences $\{ s_n \}_{n=0}^{\infty}$ rational extrapolation will give
better results than polynomial extrapolation. Let us therefore assume
that the approximant $\mathcal{S}_k (x)$ can be written as the ratio of
two polynomials of degrees $l$ and $m$, respectively:
\begin{equation}
\mathcal{S}_k (x) \; = \; \frac
{a_0 + a_1 x + a_2 x^2 + \cdots + a_l x^l}
{b_0 + b_1 x + b_2 x^2 + \cdots + b_m x^m} \; ,
\qquad k, l, m \in \mathbb{N}_0 \, .  
\label{InterpolRat}
\end{equation}
This rational function contains $l + m + 2$ coefficients $a_0,
\ldots , a_l$ and $b_0, \ldots , b_m$. However, only $l + m +1$
coefficients are independent since they are determined only up to a
common nonvanishing factor. Usually, one requires either $b_0 = 1$ or
$b_m = 1$. Consequently, the $k+1$ interpolation conditions
(\ref{IntPolCond}) will determine the coefficients $a_0, \ldots, a_l$
and $b_0, \ldots, b_m$ provided that $k = l + m$ holds.  

The extrapolation point $x_{\infty}=0$ is also the most obvious choice
in the case of rational extrapolation. Extrapolation to $x_{\infty}=0$
implies that the interpolation points $\{ x_n \}_{n=0}^{\infty}$ in
(\ref{InterpolRat}) have to satisfy (\ref{x_n2zero}) and that the
approximation to the limit is to be identified with the ratio $a_0/b_0$
of the constant terms of the polynomials in (\ref{InterpolRat}).

If $l = m$ holds in (\ref{InterpolRat}), extrapolation to infinity is
also possible. In that case the interpolation points $\{ x_n
\}_{n=0}^{\infty}$ have to satisfy
\begin{subequations}
\label{x_n2inf}
\begin{eqnarray}
& 0 < x_0 < x_1 < \cdots < x_m < x_{m+1} < \cdots \, ,  
\\
& {\displaystyle \lim_{n \to \infty} \; x_n \; = \; \infty} \, . 
\end{eqnarray}
\end{subequations}
In the case of extrapolation to infinity, only the coefficients $a_l$
and $b_l$ of the polynomials in (\ref{InterpolRat}), that are
proportional to highest power $x^l$, contribute. Consequently, the
approximation to the limit has to be identified with the ratio
$a_l/b_l$.

As in the case of polynomial interpolation, several different algorithms 
for the computation of rational interpolants are known. A discussion of
the relative merits of these algorithms as well as a survey of the
relevant literature can be found in \cite[Chapter III]{CuyWuy1987}.

The most frequently used rational extrapolation technique is probably
Wynn's rho algorithm \cite{Wy1956a}:
\begin{subequations}
\label{RhoAl}
\begin{eqnarray}
\rho_{-1}^{(n)} & \; = \; & 0 \, ,
\qquad \rho_0^{(n)} \; = \; s_n \, , \qquad n \in \mathbb{N}_0 \, , \\
\rho_{k+1}^{(n)} & \; = \; & \rho_{k-1}^{(n+1)} \, + \,
\frac { x_{n+k+1} - x_n} {\rho_{k}^{(n+1)} - \rho_{k}^{(n)}} \, ,
\qquad k, n \in \mathbb{N}_0 \, .
\end{eqnarray}
\end{subequations}     
Formally, the only difference between Wynn's epsilon algorithm
(\ref{eps_al}) and Wynn's rho algorithm is that the rho algorithm
involves a sequence of interpolation points $\{ x_n \}_{n=0}^{\infty}$
which have to satisfy (\ref{x_n2inf}). As in the case of the epsilon
algorithm, only the elements $\rho_{2 k}^{(n)}$ with even subscripts
serve as approximations to the limit. The elements $\rho_{2 k +
1}^{(n)}$ with odd subscripts are only auxiliary quantities which
diverge if the whole process converges.

Wynn's rho algorithm is designed to compute even-order convergents of
Thiele's interpolating continued fraction \cite{Thi1909} and to
extrapolate them to infinity. The even-order convergents are rational
functions of the following type:
\begin{equation}
\mathcal{S}_{2 k} (x) \; = \; \frac
{a_k x^k + a_{k-1} x^{k-1} + \cdots + a_1 x + a_0 }
{b_k x^k + b_{k-1} x^{k-1} + \cdots + b_1 x + b_0} \, ,
\qquad k \in \mathbb{N}_0 \, . 
\end{equation}
Thus, the ratio $a_k/b_k$ is to be identified with the approximation to
the limit.

The epsilon and the rho algorithm have complementary features. The
epsilon algorithm is a powerful accelerator for linear convergence and
is also able to sum many divergent alternating series, whereas the rho
algorithm fails to accelerate linear convergence and is not able to sum
divergent series. However, it is a very powerful accelerator for some
logarithmically convergent sequences.

The properties of Wynn's rho algorithm are for example discussed in
books by Brezinski (see \cite[pp.\ 102 - 106]{Bre1977} and \cite[pp.\ 96
- 102]{Bre1978} and Wimp \cite[pp.\ 168 - 169]{Wi1981}. In these books
the connection of the $\rho$ algorithm with interpolating continued
fractions is emphasized and it is also shown that the transforms
$\rho_{2 k}^{(n)}$ can be represented as the ratio of two
determinants. Moreover, there is an article by Osada \cite{Osa1996},
discussing its convergence properties. But otherwise, relatively little
seems to be known about this sequence transformation.

In the case of the rho algorithm, we can proceed as in the case of the
Aitken formula (\ref{AitFor_1}) and construct an iterated
transformation. For $k = 1$, we obtain from (\ref{RhoAl}):
\begin{equation}
\rho_2^{(n)} \; = \; s_{n+1} \, + \, \frac
{(x_{n+2} - x_n) [ \Delta s_{n+1} ] [ \Delta s_n ] }
{ [ \Delta x_{n+1} ] [ \Delta s_{n} ] -
[ \Delta x_n ] [ \Delta s_{n+1} ] }
\, , \qquad n \in \mathbb{N}_0 \, .   
\label{rho_2}
\end{equation} 
This expression, which can be viewed to be a kind of weighted $\Delta^2$
formula, can be iterated yielding  \cite[Section 6.3]{We1991}
\begin{subequations}
\label{RhoIt}
\begin{eqnarray}
{\cal W}_{0}^{(n)} & = & s_n \, ,
\qquad n \in \mathbb{N}_0 \, , 
\\
{\cal W}_{k+1}^{(n)} & = & {\cal W}_{k}^{(n+1)} \,
+ \frac
{ (x_{n + 2 k + 2} - x_n ) \bigl[ \Delta {\cal W}_k^{(n+1)}
\bigr] \bigl[ \Delta {\cal W}_k^{(n)} \bigr] }
{(x_{n + 2 k + 2} - x_{n+1}) \bigl[ \Delta {\cal W}_k^{(n)} \bigr] -
(x_{n + 2 k + 1} - x_n) \bigl[ \Delta {\cal W}_k^{(n+1)} \bigr] } \, , 
\nonumber \\
& & k, n \in \mathbb{N}_0 \, . 
\end{eqnarray} 
\end{subequations}
This is not the only possibility of iterating $\rho_2^{(n)}$. However,
in \cite{We1991} it was shown that other iterations of $\rho_2^{(n)}$ --
for instance the one derived by Bhowmick, Bhattacharya, and Roy
\cite{BhoBhaRo1989} -- are significantly less efficient than ${\cal
W}_k^{(n)}$, which has similar properties as Wynn's rho algorithm
(\ref{RhoAl}), from which it was derived.

The main problem with the practical application of those sequence
transformations, that are based upon interpolation theory, is that one
has to find a sequence $\{ x_n \}_{n=0}^{\infty}$ of interpolation
points which produces good extrapolation results for a given sequence
$\{ s_n \}_{n=0}^{\infty}$ of input data. For example, in most practical 
applications the Richardson extrapolation scheme (\ref{RichAl}) is used
in combination with the interpolation points $x_n = 1/(n+\beta)$, where
$\beta$ is a positive shift parameter. Then, $\mathcal{N}_k^{(n)}
\bigl(s_n , 1/(n+\beta) \bigr)$ possesses a closed form expression (see
for example \cite[Lemma 2.1 on p.\ 313]{MaSha1983} or \cite[Eq.\
(7.3-20)]{We1989})
\begin{eqnarray}
\lefteqn{\Lambda_k^{(n)} (\beta, s_n) \; = \; 
\mathcal{N}_k^{(n)} \bigl(s_n , 1/(n+\beta) \bigr)} \nonumber \\
\qquad & = & (-1)^k \> \sum_{j=0}^k \>
(-1)^j \> \frac {(\beta+n+j)^k} {j! \> (k-j)!} \> s_{n+j} \, ,
\qquad k, n \in \mathbb{N}_0 \, ,
\end{eqnarray}
and the recursive scheme (\ref{RichAl}) assumes the following form 
\cite[Eq.\ (7.3-21)]{We1989}
\begin{subequations}
\label{RichAlStand}
\begin{eqnarray}
\Lambda_0^{(n)} (\beta, s_n) & = & s_n \, , 
\qquad n \in \mathbb{N}_0 \, , 
\\
\Lambda_{k+1}^{(n)} (\beta, s_n) & = &
\Lambda_k^{(n+1)} (\beta, s_{n+1}) \, + \,
\frac {\beta+n} {k+1} \> \Delta \Lambda_k^{(n)} (\beta, s_n)
\, , \qquad k, n \in \mathbb{N}_0 \, .
\end{eqnarray}
\end{subequations}
Similarly, Wynn's rho algorithm (\ref{RhoAl}) and its iteration
(\ref{RhoIt}) are normally used in combination with the interpolation
points $x_n = n+1$, yielding the standard forms (see for example
\cite[Eq.\ (6.2-4)]{We1989})
\begin{subequations}
\label{RhoAlStand}
\begin{eqnarray}
\rho_{-1}^{(n)} & \; = \; & 0 \, ,
\qquad \rho_0^{(n)} \; = \; s_n \, , \qquad n \in \mathbb{N}_0 \, , \\
\rho_{k+1}^{(n)} & \; = \; & \rho_{k-1}^{(n+1)} \, + \,
\frac {k+1} {\rho_{k}^{(n+1)} - \rho_{k}^{(n)}} \, ,
\qquad k, n \in \mathbb{N}_0 \, ,
\end{eqnarray}
\end{subequations}     
and \cite[Section 6.3]{We1991}
\begin{subequations}
\label{RhoItStand}
\begin{eqnarray}
{\cal W}_{0}^{(n)} & = & s_n \, ,
\qquad n \in \mathbb{N}_0 \, , 
\\
{\cal W}_{k+1}^{(n)} & = & {\cal W}_{k}^{(n+1)} \, -
\frac
{ (2 k + 2) \bigl[ \Delta {\cal W}_k^{(n+1)} \bigr]
\bigl[ \Delta {\cal W}_k^{(n)} \bigr] }
{ (2 k + 1) \Delta^2 {\cal W}_k^{(n)} } \, ,  
\qquad k,n \in \mathbb{N}_0 \, . 
\end{eqnarray}
\end{subequations}

In the literature on convergence acceleration and extrapolation,
Richardson extrapolation as well as Wynn's rho algorithm and its
iteration are normally used in their standard forms (\ref{RichAlStand}),
(\ref{RhoAlStand}), and (\ref{RhoItStand}). However, it is important to
recognize that their general forms (\ref{RichAl}), (\ref{RhoAl}), and
(\ref{RhoIt}), respectively, which all involve unspecified interpolation
points, may be required for oligomer calculations. Let us assume that
the elements $s_n$ of a logarithmically convergent sequence can be
expressed by series expansions of the following kind,
\begin{equation}
s_n \; = \; s \, + \, (n+\beta)^{-\alpha} \, 
\sum_{j=0}^{\infty} c_j / (n+\beta)^j \, , 
\qquad n \in \mathbb{N}_0 \, ,
\label{ModSeqAlpha}
\end{equation}
where $\alpha$ is a positive decay parameter, $\beta$ is a positive
shift parameter, and the $c_j$ are unspecified coefficients. In general,
we have no \emph{a priori} reason to assume that oligomer calculations
would only produce sequences of that kind with an integral decay
parameter, and we want to accelerate the convergence of sequences of
that kind even if $\alpha$ is nonintegral.

In Theorem 14-4 of \cite{We1989}, it was shown rigorously that the
standard form (\ref{RichAlStand}) of Richardson extrapolation is able to
accelerate the convergence of the model sequence (\ref{ModSeqAlpha}) if
$\alpha = 1, 2, \ldots$, but fails if $\alpha$ is nonintegral. In the
case of Wynn's rho algorithm and its iteration, we are not aware of any
explicit proof, but there is considerable experimental evidence that the
standard forms (\ref{RhoAlStand}) and (\ref{RhoItStand}) of these
transformations also fail to accelerate the convergence of the model
sequence (\ref{ModSeqAlpha}) if the decay parameter $\alpha$ is
nonintegral (see Section 14.4 of \cite{We1989}).

If the decay parameter $\alpha$ of a sequence of the type of
(\ref{ModSeqAlpha}) is known, then Osada's variant of Wynn's rho
algorithm can be used \cite[Eq.\ (3.1)]{Os1990}:
\begin{subequations}
\label{OsRhoAl}
\begin{eqnarray}
{\bar \rho}_{-1}^{(n)} & \; = \; & 0 \, ,
\qquad {\bar \rho}_0^{(n)} \; = \; s_n \, ,
\qquad n \in \mathbb{N}_0 \, , \\
{\bar \rho}_{k+1}^{(n)} & \; = \; & {\bar \rho}_{k-1}^{(n+1)} \, + \,
\frac {k+\alpha} {{\bar \rho}_{k}^{(n+1)} - {\bar \rho}_{k}^{(n)}} \, ,
\qquad k, n \in \mathbb{N}_0 \, .
\end{eqnarray}
\end{subequations}   
For $\alpha = 1$, Osada's variant becomes identical with the standard
form (\ref{RhoAlStand}). Osada also demonstrated that his variant
accelerates the convergence of the model sequence (\ref{ModSeqAlpha}),
and that the transformation error satisfies the following asymptotic
estimate \cite[Theorem 4]{Os1990}:
\begin{equation}
{\bar \rho}_{2 k}^{(n)} \, - \, s \; = \; 
\mathrm{O} \bigl( n^{-\alpha - 2k} \bigr) \, , \qquad n \to \infty \, .  
\label{OsAsyTranErr}
\end{equation}

Osada's variant of the rho algorithm can also be iterated. From
(\ref{OsRhoAl}) we obtain the following expression for ${\bar
\rho}_2^{(n)}$ in terms of $s_n$, $s_{n+1}$, and $s_{n+2}$:
\begin{equation}
{\bar \rho}_2^{(n)} \; = \;
s_{n+1} \, - \, \frac {(\alpha + 1)} {\alpha} \,
\frac {[\Delta s_n] [\Delta s_{n+1}]} {[\Delta^2 s_n]} \, ,
\qquad n \in \mathbb{N}_0 \, . 
\end{equation}
If this iteration is done in such a way that $\alpha$ is increased by
2 with every recursive step, we obtain the following recursive scheme
\cite[Eq.\ (2.29)]{We1991} which was derived originally by Bj{\o}rstad, 
Dahlquist, and Grosse \cite[Eq.\ (2.4)]{BjDaGr1981}:
\begin{subequations}
\label{BDGal}
\begin{eqnarray}
{\overline {\cal W}}_{0}^{(n)} \; & = \; & s_n \, ,
\qquad n \in \mathbb{N}_0 \, , \\
{\overline {\cal W}}_{k+1}^{(n)} \; & = \; &
{\overline {\cal W}}_{k}^{(n+1)} \, - \,
\frac {(2 k + \alpha + 1)} {(2 k + \alpha)} \,
\frac
{\bigl[ \Delta {\overline {\cal W}}_k^{(n+1)} \bigr]
\bigl[ \Delta {\overline {\cal W}}_k^{(n)} \bigr]}
{\Delta^2 {\overline {\cal W}}_k^{(n)}} \, ,
\quad k, n \in \mathbb{N}_0 \, . \quad
\end{eqnarray}
\end{subequations}          
For $\alpha = 1$, this weighted iterated $\Delta^2$ algorithm becomes
identical with the standard form (\ref{RhoItStand}) of the iteration of
Wynn's rho algorithm.

Bj{\o}rstad, Dahlquist, and Grosse showed that ${\overline {\cal
W}}_k^{(n)}$ accelerates the convergence of the model sequence
(\ref{ModSeqAlpha}), and that the transformation error satisfies the
following asymptotic estimate \cite[Eq.\ (3.1)]{BjDaGr1981}
\begin{equation}
{\overline {\cal W}}_k^{(n)} \, - \, s \; = \; 
\mathrm{O} \bigl( n^{-\alpha - 2k} \bigr) \, , \qquad n \to \infty \, .  
\label{AsyTranErrBDG}
\end{equation}
The two sequence transformations ${\bar \rho}_{2 k}^{(n)}$ and
${\overline {\cal W}}_k^{(n)}$ require the same set $s_n$, $s_{n+1}$,
$\ldots$ , $s_{n+2k}$ of elements of the model sequence
(\ref{ModSeqAlpha}) as input data. Since the asymptotic error estimates
(\ref{OsAsyTranErr}) and (\ref{AsyTranErrBDG}) are equivalent, the two
transformations are also \emph{asymptotically} equivalent, i.e., in the
limit of large indices $n$, although they of course may produce
different results for small indices $n$. Moreover, these two
transformations are also \emph{optimal}, since the elimination of two
terms of the infinite series on the right-hand side of
(\ref{ModSeqAlpha}) requires the input of the numerical values of two
sequence elements. Thus, no sequence transformation, which only uses the
numerical values of the elements of the model sequence
(\ref{ModSeqAlpha}) as input data, could improve on the asymptotic ($n
\to \infty$) error estimates (\ref{OsAsyTranErr}) and 
(\ref{AsyTranErrBDG}).

The knowledge of the decay parameter $\alpha$ is crucial for an
application of the transformations (\ref{OsRhoAl}) and (\ref{BDGal}) to
a sequence of the type of (\ref{ModSeqAlpha}).  An approximation to
$\alpha$ can be obtained with the help of the following nonlinear
transformation:
\begin{equation}
T_n \; = \; \frac
{ [\Delta^2 s_n ] \, [\Delta^2 s_{n+1} ]}
{[\Delta s_{n+1} ] \, [\Delta^2 s_{n+1} ] \, - \,
[\Delta s_{n+2} ] \, [\Delta^2 s_{n}]} \; - \; 1 \, ,
\qquad n \in \mathbb{N}_0 \, .
\label{DecPar}
\end{equation}   
The transformation $T_n$, which was first derived in a somewhat
disguised form by Drummond \cite{Dru1976} and later rederived by
Bj{\o}rstad, Dahlquist, and Grosse \cite{BjDaGr1981}, is essentially a
weighted $\Delta^3$ method. Consequently, it is potentially a very
unstable method if the relative accuracy of the input data is low, and
stability problems can never be excluded. Bj{\o}rstad, Dahlquist, and
Grosse \cite[Eq.\ (4.1)]{BjDaGr1981} also showed that
\begin{equation}
\alpha \; = \; T_n \, + \, O (1/n^2) \,,\qquad n \to \infty \, ,
\label{T_n}
\end{equation}  
if the elements of the model sequence (\ref{ModSeqAlpha}) are used as
input data.

\setcounter{equation}{0}
\section{An Example: Extrapolation of Oligomer Calculations for Unit 
Cell Energy of Polyacetylene} 
\label{Sec:Polyac}

All-trans polyacetylene with the repeat unit $-\text{(CH=CH)}-$, which
contains 14 electrons, is the simplest \emph{quasi}-onedimensional
polymer containing conjugated double bonds. As a result of the
conjugation and the pronounced anisotropy, polyacetylene exhibits a
variety of interesting properties. Included amongst these properties is
an enormous increase in conductivity induced by doping (see for
example \cite{McGeMiMoHee,HeeMacD1980}) and a very large electric field
induced polarization \cite{KirCha1997,ChaKir2000} compared to the
corresponding non-conjugated polymer polyethylene with the repeat unit
$-\text{(CH}_2\text{-CH}_2\text{)-}$. Accordingly, there is an enormous
number of articles dealing with polyacetylene. A reasonably complete
bibliography would clearly be beyond the scope of this article.

In this article, we want to demonstrate that sequence transformations
are useful numerical tools for the extrapolation of oligomer
calculations. Consequently, we are not so much interested in obtaining
highly accurate results. Instead, we want to study the typical behavior
of quantum chemical calculations on a sequence of oligomers with
increasing size. For this purpose, it is sufficient to use relatively
small basis sets which allowed us to carry out calculations on fairly
large oligomers expeditiously.

Here we some consider calculations for the ground state energy of the
oligomers $\mathrm{H} \! - \! \text{(CH=CH)}_{\mathrm{N}}
\! - \! \mathrm{H}$. These will illustrate how new and improved results 
can be obtained but, by no means exhaust the range of
possibilities. Other \emph{quasi}-onedimensional polymers, as well as
the application to other properties such as hyperpolarizabilities and
electronic transition energies, will be discussed elsewhere
\cite{KirWen2000}. The calculated total ground state energy $E_N$ of an 
oligomer consisting of $N$ repeat units is an \emph{extensive}
quantity. This means that we first have to convert $E_N$ to an
\emph{intensive} quantity before we can do extrapolations to the
infinite chain limit. The most obvious approach consists in defining
average energies or energies per repeat unit according to
\begin{equation}
E_{N}^{\mathrm{(av)}} \; = \; E_N/N \, .
\label{DefE_av}
\end{equation}
If we compare the total energy of an oligomer $\mathrm{H} \! -
\! \text{(CH=CH)}_{\mathrm{N}} \! - \! \mathrm{H}$ with the total energy 
of a hypothetical ring molecule consisting of $N$ repeat units
$-\text{(CH=CH)}-$, and if we also assume that $N$ is so large that
angular strain can be neglected, then we see that the saturation of the
terminal free valences by hydrogen atoms should lead to an energy
contribution which is nearly constant, i.e., which depends on $N$ only
quite weakly. Consequently, we may reasonably assume that the average
energy $E_{N}^{\mathrm{(av)}}$ will contain a contribution which is
proportional to $1/N$ arising from the chain ends and, therefore, can be
expressed as
\begin{equation}
E_{N}^{\mathrm{(av)}} \; = \; E_{\infty}^{\mathrm{(av)}} \, + \, C/N 
\, + \, \Phi_{N}^{\mathrm{(av)}} \, ,
\label{E_avMod}
\end{equation} 
where $C$ is an $N$-independent constant and $\Phi_{N}^{\mathrm{(av)}}$ is 
an unknown function of $N$ which should decay faster than $1/N$ as $N$
becomes large.

If we want to extrapolate the average energies successfully to the
infinite chain limit, then we must have an idea of how the unknown
function $\Phi_{N}^{\mathrm{(av)}}$ depends upon $N$. It is a relatively
obvious idea to express $\Phi_{N}^{\mathrm{(av)}}$ or equivalently
$E_{N}^{\mathrm{(av)}}$ by a power series in $1/N$ according to
\begin{equation}
E_{N}^{\mathrm{(av)}} \; = \; E_{\infty}^{\mathrm{(av)}} \, + \,
\sum_{j=1}^{\infty} \, C_j / N^j \, .
\label{PowSerE_N}
\end{equation}
If this is case, then the sequence of the average energies converges
logarithmically to the infinite chain limit according to
(\ref{LinLogConv}). As discussed in Section \ref{Sec:RichRho}, the
acceleration of logarithmic convergence is normally more difficult than
the acceleration of linear convergence. Nevertheless, the standard forms
(\ref{RichAlStand}), (\ref{RhoAlStand}), and (\ref{RhoItStand}) of
Richardson extrapolation, Wynn's rho algorithm and its iteration should
be effective accelerators for sequences of the kind of
(\ref{PowSerE_N}). 

However, there are still some unresolved theoretical
questions. According to Cioslowski and Lepetit \cite[p.\
3544]{CioLep1991}, who performed a perturbation theoretic analysis of
the $N$-dependence of Hartree-Fock oligomer energies, the power series
expansion (\ref{PowSerE_N}) for the average energy
$E_{N}^{\mathrm{(av)}}$ can not only contain powers of $1/N$ but also
nonanalytic terms of the type $\exp (-\gamma N)$ with $\gamma > 0$. For
sufficiently large values of $N$, these nonanalytic contributions will
be negligible. However, they could well have a negative effect on the
performance of convergence acceleration processes. Sequence
transformations designed for the acceleration of logarithmic convergence
are usually not good at eliminating nonanalytic contributions typical of
linear convergence, and vice versa.

Alternatively, Cui, Kertesz, and Jiang \cite{CuiKerJia1990} proposed to
compute the energy per repeat unit not via the average energies
(\ref{DefE_av}) but via the energy differences
\begin{equation}
E_{N}^{\mathrm{(dif)}} \; = \; E_{N+1} \, - \, E_N 
\; = \; (N+1) \, E_{N+1}^{\mathrm{(av)}} \, - \,
N \, E_{N}^{\mathrm{(av)}} \, .
\label{DefE_dif}
\end{equation}
The argument of Cui, Kertesz, and Jiang was that the energies of two
oligomers containing $N+1$ or $N$ monomer units, respectively, should
essentially differ by the energy of an \emph{inner} monomer unit
\cite{Kir1992}. Equivalently, this means that the effects due to the
finite size of the oligomers as well as the $1/N$ errors due to the
saturation of free terminal valences should depend on $N$ only quite
weakly and are largely canceled by forming differences.

If we now combine (\ref{E_avMod}) with (\ref{DefE_dif}), we obtain
\begin{equation}
E_{N}^{\mathrm{(dif)}} \; = \; E_{\infty}^{\mathrm{(av)}} \, + \,
\Delta \bigl[ N \, \Phi_{N}^{\mathrm{(av)}} \bigr] \, ,
\label{E_difMod}
\end{equation}
where $\Phi_{N}^{\mathrm{(av)}}$ is the function defined in
(\ref{E_avMod}) that describes the higher order effects of the
$N$-dependence of the average energies. Comparing (\ref{E_avMod}) and
(\ref{E_difMod}) and taking into account that $\Phi_{N}^{\mathrm{(av)}}$
has been assumed to decay faster than $1/N$, we see that
$E_{N}^{\mathrm{(dif)}}$ does not contain a $1/N$ contribution and
converges more rapidly than $E_{N}^{\mathrm{(av)}}$ to the infinite
chain limit. If, for example, it is assumed that
$\Phi_{N}^{\mathrm{(av)}} \sim 1/N^2$ as $N \to \infty$, we obtain
$\Delta \bigl[ N \Phi_{N}^{\mathrm{(av)}} \bigr] \sim 1/N^2$ which shows
that the energy differences converge faster than the average energies.

In order to illustrate some of the above procedures we carried out
Hartree-Fock ground state energy calculations for the oligomers
$\mathrm{H} \! - \! \text{(CH=CH)}_{\mathrm{N}} \! - \!
\mathrm{H}$ with $N \le 16$ using an \texttt{STO-3G} basis set. A fixed
geometry was used for all oligomers with ${R}_{\mathrm{C=C}} =
1.366~{\AA}$, ${R}_{\mathrm{C-C}} = 1.450~{\AA}$, ${R}_{\mathrm{C-H}} =
1.085~{\AA}$, ${\angle}_{\mathrm{C=C-C}} = 123.90~{\deg}$,
${\angle}_{\mathrm{H-C-H}} = 117.13~{\deg}$, and
${\angle}_{\mathrm{H-C-H}} = 117.13~{\deg}$. This geometry implies
${\angle}_{\mathrm{H-C=C}} = 118.97~{\deg}$, except for the terminal
hydrogen on either end, for which ${\angle}_{\mathrm{H-C=C}} =
123.90~{\deg}$. In Table I, we list the total energies $E_N$, the
average energies $E_{N}^{\mathrm{(av)}}$ defined in (\ref{DefE_av}), and
the energy differences $E_{N}^{\mathrm{(dif)}}$ defined in
(\ref{DefE_dif}). The data displayed in Table I show that the energy
differences converge indeed much more rapidly than the average energies.

\begin{table}[htb] 
\begin{center}
\begin{tabular*}{10.5cm}{l@{\extracolsep{0.5cm}}rrr} \\
\multicolumn{4}{l}{Table I: Total Hartree-Fock energies $E_N$,
average energies $E_{N}^{\mathrm{(av)}}$,} \\
\multicolumn{4}{l}{and energy differences $E_N^{\mathrm{(dif)}}$ using 
an \texttt{STO-3G} basis set.} \\
\multicolumn{4}{l}{For geometry see the text.} \\
[1\jot] \hline \hline \rule{0pt}{4\jot}%
$N$ & \multicolumn{1}{c}{$E_{N}$} 
& \multicolumn{1}{c}{$E_{N}^{\mathrm{(av)}}$}
& \multicolumn{1}{c}{$E_{N}^{\mathrm{(dif)}}$}
\\ [1\jot] \hline \rule{0pt}{4\jot}%
1  & $-77.0672438490$ & $-77.0672438490$ & $-75.943944441$ \\	      
2  & $-153.011188290$ & $-76.5055941450$ & $-75.945112888$ \\ 
3  & $-228.956301178$ & $-76.3187670593$ & $-75.945528271$ \\ 
4  & $-304.901829449$ & $-76.2254573623$ & $-75.945641947$ \\ 
5  & $-380.847471396$ & $-76.1694942792$ & $-75.945676982$ \\ 
6  & $-456.793148378$ & $-76.1321913963$ & $-75.945688518$ \\ 
7  & $-532.738836896$ & $-76.1055481280$ & $-75.945692475$ \\ 
8  & $-608.684529371$ & $-76.0855661714$ & $-75.945693869$ \\ 
9  & $-684.630223240$ & $-76.0700248044$ & $-75.945694368$ \\ 
10 & $-760.575917608$ & $-76.0575917608$ & $-75.945694549$ \\ 
11 & $-836.521612157$ & $-76.0474192870$ & $-75.945694615$ \\ 
12 & $-912.467306772$ & $-76.0389422310$ & $-75.945694639$ \\ 
13 & $-988.413001411$ & $-76.0317693393$ & $-75.945694649$ \\ 
14 & $-1064.35869606$ & $-76.0256211471$ & $-75.945694650$ \\ 
15 & $-1140.30439071$ & $-76.0202927140$ & $-75.945694650$ \\ 
16 & $-1216.25008536$ & $-76.0156303350$ &                 \\	
[1\jot] \hline \hline \rule{0pt}{4\jot}%
\end{tabular*}
\end{center}
\end{table}

In the next step, we want to analyze whether the average energies and
the energy differences behave like the elements of the model sequence
(\ref{ModSeqAlpha}), i.e., whether we can identify a well-defined decay
parameter $\alpha$. For that purpose, we use the transformation $T_n$
defined in (\ref{DecPar}). All sequence transformations discussed in
this article are gauged in such a way that the indices of the elements
of the sequence to be transformed start with zero. Thus, we use the
average energies and the energy differences as input data for $T_n$
according to either $s_n = E_{n+1}^{\mathrm{(av)}}$ with $0 \le n \le
16$ or $s_n = E_{n+1}^{\mathrm{(dif)}}$ with $0 \le n \le 14$.

\begin{table}[htb] 
\begin{center}
\begin{tabular*}{11.2cm}{l@{\extracolsep{0.3cm}}rrrr} \\
\multicolumn{5}{l}{Table II: Approximations to the
decay parameter $\alpha$ of the average} \\
\multicolumn{5}{l}{energies $E_{n+1}^{\mathrm{(av)}}$ and the energy 
differences $E_{n+1}^{\mathrm{(dif)}}$ according to 
(\protect\ref{DecPar}).} \\
[1\jot] \hline \hline \rule{0pt}{4\jot}%
$n$  
& \multicolumn{1}{c}{$E_{n+1}^{\mathrm{(av)}}$} 
& \multicolumn{1}{c}{$T_n$} 
& \multicolumn{1}{c}{$E_{n+1}^{\mathrm{(dif)}}$}
& \multicolumn{1}{c}{$T_n$} 
\\ [1\jot] \hline \rule{0pt}{4\jot}%
0  & $-77.0672438490$ & $1.0026524$ & $-75.943944441$ & $-6.7203517$ \\	      
1  & $-76.5055941450$ & $0.9972079$ & $-75.945112888$ & $ 13.549818$ \\ 
2  & $-76.3187670593$ & $0.9976702$ & $-75.945528271$ & $ 21.022075$ \\ 
3  & $-76.2254573623$ & $0.9984106$ & $-75.945641947$ & $ 31.065636$ \\ 
4  & $-76.1694942792$ & $0.9990241$ & $-75.945676982$ & $ 44.885592$ \\ 
5  & $-76.1321913963$ & $0.9994399$ & $-75.945688518$ & $ 72.270674$ \\ 
6  & $-76.1055481280$ & $0.9996933$ & $-75.945692475$ & $ 84.907033$ \\ 
7  & $-76.0855661714$ & $0.9998391$ & $-75.945693869$ & $ 210.38728$ \\ 
8  & $-76.0700248044$ & $0.9999177$ & $-75.945694368$ & $-403.50000$ \\ 
9  & $-76.0575917608$ & $0.9999589$ & $-75.945694549$ & $ 6.0000000$ \\ 
10 & $-76.0474192870$ & $0.9999827$ & $-75.945694615$ & $-2.6578947$ \\ 
11 & $-76.0389422310$ & $0.9999829$ & $-75.945694639$ & $-10.000000$ \\ 
12 & $-76.0317693393$ & $0.9999976$ & $-75.945694649$ &              \\ 
13 & $-76.0256211471$ &             & $-75.945694650$ &              \\ 
14 & $-76.0202927140$ &             & $-75.945694650$ &              \\ 
15 & $-76.0156303350$ &             &                 &              \\	
[1\jot] \hline \hline \rule{0pt}{4\jot}%
\end{tabular*}
\end{center}
\end{table}

The results in the third column of Table II show quite clearly that the
average energies possess a contribution which is proportional to
$1/N$. On the basis of these results alone, one might be tempted to
generalize this observation and conclude that the average energies can
be expressed by expansions in powers of $1/N$ according to
(\ref{PowSerE_N}). However, the approximations to the decay parameter of
the energy differences in the fifth column show an erratic behavior and
there is no indication that they might converge. Thus, we have no reason
to assume that energy differences decay like a fixed power of $1/N$.
In view of {\ref{DefE_dif}) and (\ref{E_difMod}) this observation
implies that the unknown function $\Phi_{N}^{\mathrm{(av)}}$ defined in
(\ref{E_avMod}), which describes the higher order effects of the
truncation errors in the average energies, also does not decay like a
fixed power of $1/N$. Consequently, our observation contradicts the
perturbation theoretic analysis of Cioslowski and Lepetit
\cite{CioLep1991} who concluded that the average energies can be
expressed by a power series in $1/N$ according to (\ref{PowSerE_N})
plus possibly some nonanalytic contributions.

If we want to apply sequence transformations to speed up the convergence
of oligomer calculations to the infinite chain limit, we have to assume
that the $N$-dependence of the truncation errors in the energy
differences possesses a structure which is sufficiently regular to be
detected and utilized in convergence acceleration processes. The energy
differences in the last column of Table I showed a remarkably rapid
convergence, indicating that the truncation errors might converge
exponentially. A rigorous theoretical proof of this assumption is of
course not in sight. So, we have to provide numerical evidence for our
conjecture.

There is a simple ratio test related to the ratio test for infinite
series, which can provide evidence for exponential decay. If we form the
ratio
\begin{equation}
\mathcal{R}_n \; = \; 
\frac{s_{n+2} - s_{n+1}}{s_{n+1} - s_{n}} \; = \; 
\frac{\Delta s_{n+1}}{\Delta s_{n}}
\label{RatioTest}
\end{equation}
of the elements of the model sequence (\ref{AitModSeq}), then we obtain
$\mathcal{R}_n = \lambda$. If the input data do not have the simple form
of the elements of this model sequence, but nevertheless decay
exponentially, then we only get a constant as $n \to \infty$. In Table
III, we apply the ratio test (\ref{RatioTest}) to the input data $s_n =
E_{n+1}^{\mathrm{(dif)}}$ with $0 \le n \le 11$. The results in the last
column of Table III seem to converge, albeit slowly. Therefore, they
provide evidence that the energy differences indeed converge
exponentially to the infinite chain limit.

\begin{table}[htb] 
\begin{center}
\begin{tabular*}{8.2cm}{l@{\extracolsep{1.9cm}}rr} \\
\multicolumn{3}{l}{Table III: Ratio test (\protect\ref{RatioTest}) of
the energy differences} \\
[1\jot] \hline \hline \rule{0pt}{4\jot}%
$n$ 
& \multicolumn{1}{c}{$E_{n+1}^{\mathrm{(dif)}}$} 
& \multicolumn{1}{c}{$\mathcal{R}_n$} 
\\ [1\jot] \hline \rule{0pt}{4\jot}%
0 & $-75.943944441$ & $0.3555$ \\	      
1 & $-75.945112888$ & $0.2737$ \\ 
2 & $-75.945528271$ & $0.3082$ \\ 
3 & $-75.945641947$ & $0.3293$ \\ 
4 & $-75.945676982$ & $0.3430$ \\ 
5 & $-75.945688518$ & $0.3523$ \\ 
6 & $-75.945692475$ & $0.3580$ \\ 
7 & $-75.945693869$ & $0.3627$ \\ 
8 & $-75.945694368$ & $0.3646$ \\ 
9 & $-75.945694549$ & $0.3636$ \\ 
10& $-75.945694615$ &          \\ 
11& $-75.945694639$ &          \\ 
[1\jot] \hline \hline \rule{0pt}{4\jot}%
\end{tabular*}
\end{center}
\end{table}

There is another approach which can provide evidence that the truncation
errors of the energy differences indeed decay exponentially. If we apply
sequence transformations, that are known to accelerate linear
convergence but fail in the case of logarithmic convergence like
Aitken's iterated $\Delta^2$ process (\ref{It_Aitken}) or Wynn's epsilon
algorithm (\ref{eps_al}), and observe an acceleration of convergence,
then we have very strong evidence for our conjecture.

There is a problem with the acceleration of the convergence of the
energy differences. The last column of Table I shows that the energy
differences $E_{N}^{\mathrm{(dif)}}$ have for $N = 14$ already converged
to all digits given. Thus, there are no additional digits left which
could provide information to the extrapolation methods. All sequence
transformations form higher weighted differences of the input data, and
so does Aitken's iterated $\Delta^2$ process or Wynn's epsilon
algorithm. If, however, the input data have almost reached convergence
and there are no additional digits left, numerical instabilities or even
completely nonsensical results can never be excluded.

If the input data $s_{0}$, $s_{1}$, \ldots, $s_{2k}$ are known, then the
transform with the highest possible subscript, that can be computed with
the help Wynn's epsilon algorithm (\ref{eps_al}), is
$\epsilon_{2k}^{(0)}$. In general, this is the most efficient approach
and one can normally expect the best results. If, however, the input
data behave like the energy differences and differ at most with respect
to a few digits, it is recommended that one apply the epsilon algorithm
differently. Instead, we start from the initial string
$\epsilon_{n}^{(2)} = E_{n+1}^{\mathrm{(dif)}}$ with $0
\le n \le n_{\mathrm{max}}$ and compute the string $\epsilon_{n}^{(2)}$
with $0 \le n \le n_{\mathrm{max}}-2$. Then we analyze the convergence
and reliability of this string and compute, if necessary and possible,
the next string $\epsilon_{n}^{(4)}$ with $0 \le n \le
n_{\mathrm{max}}-4$. This process can be repeated until either
convergence has been achieved or the accumulation of rounding errors
makes it impossible to proceed.

\begin{table}[htb] 
\begin{center}
\begin{tabular*}{13cm}{l@{\extracolsep{0.45cm}}rrrr} \\
\multicolumn{5}{l}{Table IV: Epsilon extrapolation of
the energy differences.} \\
[1\jot] \hline \hline \rule{0pt}{4\jot}%
$n$ 
& \multicolumn{1}{c}{$\epsilon_{0}^{(n)} = E_{n+1}^{\mathrm{(dif)}}$} 
& \multicolumn{1}{c}{$\epsilon_{2}^{(n)}$} 
& \multicolumn{1}{c}{$\epsilon_{4}^{(n)}$} 
& \multicolumn{1}{c}{$\epsilon_{6}^{(n)}$} 
\\ [1\jot] \hline \rule{0pt}{4\jot}%
0  & $-75.943944441$ & $-75.945757392$ & $-75.945691527$ & $-75.945694631$ \\  
1  & $-75.945112888$ & $-75.945684777$ & $-75.945694512$ & $-75.945694655$ \\
2  & $-75.945528271$ & $-75.945692590$ & $-75.945694634$ & $-75.945694651$ \\
3  & $-75.945641947$ & $-75.945694181$ & $-75.945694652$ & $-75.945694652$ \\
4  & $-75.945676982$ & $-75.945694541$ & $-75.945694651$ & $-75.945694653$ \\
5  & $-75.945688518$ & $-75.945694627$ & $-75.945694654$ & $-75.945694652$ \\
6  & $-75.945692475$ & $-75.945694646$ & $-75.945694653$ & $-75.945694653$ \\
7  & $-75.945693869$ & $-75.945694652$ & $-75.945694653$ & $-75.945694763$ \\
8  & $-75.945694368$ & $-75.945694653$ & $-75.945694653$ &		   \\
9  & $-75.945694549$ & $-75.945694653$ & $-75.945694654$ &		   \\
10 & $-75.945694615$ & $-75.945694656$ &                 &		   \\
11 & $-75.945694639$ & $-75.945694650$ &                 &		   \\
12 & $-75.945694649$ & 		       &                 &		   \\
13 & $-75.945694650$ &                 &                 & 		   \\
[1\jot] \hline \hline \rule{0pt}{4\jot}%
\end{tabular*}
\end{center}
\end{table}

In Table IV, we show the effect of Wynn's epsilon algorithm
(\ref{eps_al}) on the energy differences. The extrapolation results
provide strong evidence that the truncation errors of the energy
differences converge exponentially. For example,
$\epsilon_{6}^{(1)} = -75.945694655$, which requires the initial values
$\epsilon_{0}^{(1)} = E_{2}^{\mathrm{(dif)}}$, $\epsilon_{0}^{(2)} =
E_{3}^{\mathrm{(dif)}}$, \ldots, $\epsilon_{0}^{(7)} =
E_{8}^{\mathrm{(dif)}}$ for its computation, reproduces all but the last
digit of $\epsilon_{0}^{(13)} = E_{14}^{\mathrm{(dif)}} =
-75.945694650$. Thus, the epsilon algorithm makes it unnecessary to
compute the oligomer energies $E_N$ with $9 \le N \le 14$, which is
indeed a remarkable achievement.

We also used Aitken's iterated $\Delta^2$ process (\ref{It_Aitken}) for
the extrapolation of the energy differences. These extrapolation
results were virtually identical with the epsilon extrapolation results
presented in Table IV.

When we performed oligomer calculations for polyacetylene using more
extensive basis sets and other geometries, or for other
\emph{quasi}-onedimensional polymers, we again observed an exponential
decay of the energy differences \cite{KirWen2000}. Thus, this
exponential decay seems to be a general feature which is not restricted
to polyacetylene.

\setcounter{equation}{0}
\section{Summary and Conclusions} 
\label{Sec:SumConclu}

\emph{Quasi}-onedimensional stereoregular polymers like
polyacetylene are currently of considerable interest, not only
scientifically but also because of numerous possible technological
applications. There are basically two different approaches for doing
\emph{ab initio} electronic structure calculations on such systems 
which have largely complementary advantages and disadvantages. One
method, the so-called crystal orbital method, is based on the concepts
of solid state theory. It uses periodic boundary conditions and leads to
a band structure treatment. Accordingly, this approach is well suited
to describe those features that depend crucially on the in principle
infinite extension of the polymer. However, there are still some open
computational problems with the crystal orbital method. For example, it
is very difficult to properly take into account the effect of
nonperiodic perturbations such as the interaction with an electric
field or to incorporate high level electron correlation effects.

The other method is essentially a quantum chemical approach. It
approximates the polymer by oligomers consisting of a finite number of
monomer units, i.e., by molecules of finite size. In this way, the
highly developed technology of quantum chemical molecular \emph{ab
initio} programs can be used. Thus, electron correlation can be treated
using any of the techniques available in these programs. In addition,
nonlinear optical properties may be calculated via this approach. At the
present time, this cannot be done reliably within the framework of
crystal orbital theory even at the Hartree-Fock level.

Unfortunately, oligomers of finite size are not necessarily able to
model those features of a polymer which crucially depend on its in
principle infinite extension since finite oligomer calculations can
suffer badly from slow convergence to the infinite chain limit.
However, one can perform electronic structure calculations for a
sequence of oligomers with an increasing number of monomer units and try
to determine the limit of this sequence with the help of suitable
extrapolation methods.

Many different approaches for doing the extrapolations are possible. In
this article, we use sequence transformations. The mechanism, whereby
sequence transformations accomplish an acceleration of convergence, can
be understood via the Euler-Maclaurin formula which is able to produce
rapidly convergent approximations to the truncation errors of some
slowly convergent infinite series like the Dirichlet series for the
Riemann zeta function. Although the Euler-MacLaurin formula is extremely
powerful in some cases, it is an analytic convergence acceleration
method and cannot be used if only the numerical values of some elements
of a slowly convergent sequence are known. Sequence transformations
also try to accomplish an acceleration by constructing approximations to
the truncation errors of the elements of a slowly convergent sequence
which are subsequently eliminated from the input data. However, they do
this by by purely numerical means.

Two practically very important classes of sequences are linearly and
logarithmically convergent sequences. The partial sums of a power series
with a nonzero, but finite radius of convergence constitute a typical
example of a linearly convergent sequence, whereas the partial sums of
the Dirichlet series for the Riemann zeta function, which is notorious
for slow convergence, is a typical example of a logarithmically
convergent sequence. Thus, the properties of these two classes of slowly
convergent sequences differ considerably, and so do the sequence
transformations which are able to accelerate effectively either linear
or logarithmic convergence. In general, the acceleration of logarithmic
convergence is more difficult than the acceleration of linear
convergence, both theoretically and practically.

Two powerful accelerators for linear convergence are Aitken's iterated
$\Delta^2$ process and Wynn's epsilon algorithm which, in the case of
the partial sums of a power series, produces Pad\'e
approximants. However, these two transformations fail to accelerate
logarithmic convergence, whereas the powerful accelerators for
logarithmic convergence -- for example Richardson extrapolation or
Wynn's rho algorithm and its iteration -- fail to accelerate linear
convergence.

There are also some other problems with the acceleration of logarithmic
convergence. In the literature on convergence acceleration, the
transformations mentioned above are usually treated in their so-called
standard forms. However, the standard forms are not able to accelerate
the convergence of those logarithmically convergent sequences whose
elements decay like a negative nonintegral power of the index. There is
no \emph{a priori} reason why sequences with a nonintegral decay
parameter cannot occur in the context of oligomer
calculations. Therefore, we also discuss transformations like Osada's
variant of Wynn's rho algorithm or the weighted $\Delta^2$ algorithm
introduced by Bj{\o}rstad, Dahlquist, and Grosse, which are both able to
accelerate the convergence of sequences of that kind.

We carried out calculations of the Hartree-Fock ground state energy for
oligomers $\mathrm{H} \! - \! \text{(CH=CH)}_{\mathrm{N}} \! - \! 
\mathrm{H}$ using a \texttt{STO-3G} basis set with a fixed geometry, and
extrapolated these energies to the infinite chain limit corresponding to
all-trans polyacetylene. The extrapolations were greatly facilitated by
the fact that no accelerators were needed for the more troublesome
logarithmic convergence. Contrary to previous belief we showed that
the truncation error for the average energy per monomer unit contains
linear term in $1/N$, where $N$ is the number of monomer units, but no
higher powers. This leading contribution can be eliminated easily by
forming differences. and the higher order contributions converge
exponentially. Thus, accelerators for linear convergence sufficed for
our purposes and Wynn's epsilon algorithm proved very effective in this
regard. The ground state energy is just one property of interest. Future
applications will deal with nonlinear optical properties and electronic
transition energies which have a different dependence upon N and are
more problematic from the practical point of view since they converge
more slowly.

\section*{Acknowledgments}

EJW thanks the Fonds der Chemischen Industrie for financial support.

\renewcommand{\baselinestretch}{1}

\end{document}